\documentclass{article}
\usepackage[T2A]{fontenc}
\usepackage[cp1251]{inputenc}
\usepackage[english,russian]{babel}
\usepackage[tbtags]{amsmath}
\usepackage{amsfonts,amssymb,mathrsfs,amscd}
\numberwithin{equation}{section}
\newtheorem{remark}{Remark}[section]
\newtheorem{lemma}{Lemma}[section]
\newtheorem{theorem}{Theorem}[section]

\newtheorem{conclusion}{Conclusion}

\newcommand{\col}{\mathop{\rm col}}

\begin{document}
\begin{Large}
\thispagestyle{empty}
\begin{center}
{\bf Inverse scattering problem for a third-order differential operator with double potential\\
\vspace{5mm}
V. A. Zolotarev}\\

B. Verkin Institute for Low Temperature Physics and Engineering
of the National Academy of Sciences of Ukraine\\
47 Nauky Ave., Kharkiv, 61103, Ukraine

\end{center}
\vspace{5mm}

{\small {\bf Abstract.} Direct and inverse scattering problems for a third-order self-adjoint differential operator on the whole axis are studied. This operator is the sum of three summands: operator of third derivative, operator of multiplication by a function, and operator of multiplication by derivative of a function. For the solution of the inverse scattering problem, two closed systems of linear integral equations are obtained. Knowing solutions to these systems, using explicit formulas, methods of restoration of both potentials on half-axes $\mathbb{R}_\pm$ are specified.}
\vspace{5mm}


\vspace{5mm}

\begin{center}
{\bf Introduction}
\end{center}
\vspace{5mm}

Inverse scattering problem for second-order differential operators (Sturm -- Liouville problem) owes its development to classical works by V. A. Marchenko \cite{1} and M. G. Krein \cite{2}. General concepts of this scattering problem and its further development are given in \cite{3} -- \cite{5}. In the 1970s, an effective and unexpected application of the method of inverse scattering problem for integration of non-linear equations in partial derivatives was discovered. Search for L-A Lax pairs for many nonlinear equations led to a second-order operator $L$ (Sturm-Liouville). However, for the Degasperis -- Processi equation, the operator $L$ is of third order \cite{7}, \cite{8}. Inverse scattering problem for a third-order operator which is equal to the sum of the operator of the third derivative and of the operator of multiplication by a function is solved in the works \cite{9,10}.

In works \cite{11} -- \cite{16} inverse problems (spectral and scattering) corresponding to a third-order equation of the form
\begin{equation}
y'''(x)+p(x)y'(x)+q(x)y(x)=\lambda y(x)\label{eq0.1}
\end{equation}
are studied. Presence of two potentials $p(x)$ and $q(x)$ caused certain difficulties.

This work solves inverse scattering problem on the whole axis of the form
\begin{equation}
-iy'''(x)+i(p(x)y(x))'+ip(x)y'(x)+q(x)y(x)=\lambda^3y(x)\quad(x\in\mathbb{R});\label{eq0.2}
\end{equation}
which corresponds to the formally self-adjoint operator
\begin{equation}
L_{p,q}=\left(i\frac d{dx}\right)^3+i\left(\frac d{dx}p(x)+p(x)\frac d{dx}\right)+q(x)\label{eq0.3}
\end{equation}
where $p(x)$ and $q(x)$ are real functions, besides $p(x)$ is differentiable and
\begin{equation}
\int\limits_{\mathbb{R}}|q(x)|^2e^{2a|x|}dx<\infty;\quad\int\limits_\mathbb{R}|p(x)|^2e^{2a|x|}dx<\infty;\quad\int\limits_\mathbb{R}|p'(x)|^2e^{2a|x|}dx<\infty.\label{eq0.4}
\end{equation}
The manuscript consists of three sections. Section 1 proves the existence of Jost solutions for equation \eqref{eq0.2} with proper asymptotics at ``$+\infty$'' and at ``$-\infty$''. Analytic (relative to $\lambda$) properties of such Jost solutions are studied and their behavior for $\lambda\rightarrow\infty$ inside corresponding sectors is described.

Section 2 solves the inverse scattering problem for waves incident from ``$+\infty$''. Main properties of the transition matrix are studied. Riemann boundary value problems for holomorphic functions in the sectors with jumps on the rays separating these sectors are obtained. The single Riemann boundary value problem at three sectors with angle $2\pi/3$ is synthesized from these boundary value problems. To study the bound states, zeros of a holomorphic at a sector function are studied and it is proved that there is a finite number of zeros and each zero is of multiplicity $2$, besides, cubes of these zeros are eigenvalues of the operator $L_{p,q}$ \eqref{eq0.3}. This section obtains a closed system of linear singular integral equations, which, in this case, is an analogue of a well-known Marchenko equation. Scattering coefficients and points corresponding to bound states are independent parameters of this system. Knowing the solution to this system of singular equations, using simple formulas, both potentials $p(x)$ and $q(x)$ are restored at the right half-axis.

Section 3 deals with the inverse scattering problem for waves incident from ``$-\infty$''. This scattering problem is dual to the scattering problem studied in Section 2. This section also constructs a closed system of linear singular integral equations (scattering coefficients and points corresponding to bound states are independent parameters of this system). As in Section 2, potentials $p(x)$ and $q(x)$  on the left half-axis $x\in\mathbb{R}$ are restored from the solution to this equation system, using basic formulas.

In conclusion, note that methods of solution of both scattering problems are based upon constructions and results of works \cite{9,10}.

\section{Preliminary information. Jost solutions}\label{s1}

{\bf 1.1} The third-order differential equation
\begin{equation}
-iy'''(\lambda,x)+i[(p(x)y(\lambda,x))'+p(x)y'(\lambda,x)]+q(x)y(\lambda,x)=\lambda^3y(\lambda,x)\label{eq1.1}
\end{equation}
($x\in\mathbb{R}$, $\lambda\in\mathbb{C}$) generates the self-adjoint operator
\begin{equation}
L_{p,q}=-i\frac{d^3}{dx^3}+i\left[\frac d{dx}p(x)+p(x)\frac d{dx}\right]+q(x)\label{eq1.2}
\end{equation}
in the space $L^2(\mathbb{R})$ with domain $\mathfrak{D}(L_{p,q})=W_2^3(\mathbb{R})$. Functions $p(x)$ and $q(x)$ are real and $p(x)\in\mathbb{C}^1(\mathbb{R})$, besides,
\begin{equation}
\int\limits_{\mathbb{R}}|q(x)|^2e^{2a|x|}<\infty;\quad\int\limits_{\mathbb{R}}|p(x)|^2e^{2a|x|}dx<\infty;\quad\int\limits_{\mathbb{R}}|p'(x)|^2e^{2a|x|}dx<\infty\label{eq1.3}
\end{equation}
($a>0$).

\begin{remark}\label{r1.1}
Conditions \eqref{eq1.3} imply that functions $|q(x)|e^{b|x|}$, $|p(x)|e^{b|x|}$, $|p'(x)|e^{b|x|}$ belong to $L^1(\mathbb{R})\cap L^2(\mathbb{R})$ for all $b$ such that $b<a$.
\end{remark}

Relation \eqref{eq1.3} implies that equation \eqref{eq1.1}, for $|x|\rightarrow\infty$, becomes
\begin{equation}
-iy'''(\lambda,x)=\lambda^3y(x)\label{eq1.4}
\end{equation}
which has three linearly independent solutions $\{e^{-i\lambda\zeta_kx}\}_0^2$ where $\{\zeta_k\}_0^2$ are roots of the cubic equation $z^3=1$:
\begin{equation}
\zeta_0=1,\quad\zeta_1=\left(-\frac12+i\frac{\sqrt3}2\right),\quad\zeta_2=\left(-\frac12-i\frac{\sqrt3}2\right).\label{eq1.5}
\end{equation}
An important role is played by another system of fundamental solutions to equation \eqref{eq1.4} \cite{9,10}:
\begin{equation}
s_p(z)=\frac13\sum\limits_k\frac1{\zeta_k^p}e^{z\zeta_k}\quad(0\leq p\leq2)\label{eq1.6}
\end{equation}
which is analogous to classic trigonometric functions and plays a key role in the solution of inverse problems for third-order equations \cite{9,10}. The main properties of the functions $\{s_p(z)\}_0^2$ \eqref{eq1.6} are given in \cite{9,10}.

By $L_{\zeta_k}$, we denote straight lines in $\mathbb{C}$ in direction of unit vectors $\{\zeta_k\}_0^2$ \eqref{eq1.5}
\begin{equation}
L_{\zeta_k}\stackrel{\rm def}{=}\{x\zeta_k:x\in\mathbb{R}\}\quad(0\leq k\leq2),\label{eq1.7}
\end{equation}
and let $l_{\zeta_k}$ be the rays in direction of $\zeta_k$ from the origin, and $\widehat{l}_{\zeta_k}$, the rays in direction of $\zeta_k$ coming to the origin,
\begin{equation}
l_{\zeta_k}\stackrel{\rm def}{=}\{x\zeta_k:x\in\mathbb{R}_+\};\quad\widehat{l}_{\zeta_k}=L_{\zeta_k}\backslash l_{\zeta_k}\stackrel{\rm def}{=}\{x\zeta_k:x\in\mathbb{R}_-\}.\label{eq1.8}
\end{equation}
Straight lines $\{L_{\zeta_k}\}_0^2$ \eqref{eq1.7} divide plane $\mathbb{C}$ into six sectors,
\begin{equation}
S_p\stackrel{\rm def}{=}\left\{z\in\mathbb{C}:\frac{2\pi}6p<\arg z<\frac{2\pi}6(p+1)\right\}\quad(0\leq p\leq5).\label{eq1.9}
\end{equation}

Consider the Cauchy problem
\begin{equation}
-iy'''(x)=\lambda^3y(x)+f(x);\quad y(0)=y_0;\,y'(0)=y_1,\,y''(0)=y_2\quad(f\in L^2(\mathbb{R})).\label{eq1.10}
\end{equation}
It is easy to see [ , ] that solution to the homogeneous problem \eqref{eq1.10} ($f\equiv0$) is
$$y_0(\lambda,x)=y_0s_0(-i\lambda x)+y_1\frac{s_1(-i\lambda x)}{(-i\lambda)}+y_2\frac{s_2(-i\lambda x)}{(-i\lambda)^2}.$$
Hence, by the method of variation of arbitrary constants, we find the solution to nonhomogeneous ($f\not\equiv0$) Cauchy problem \eqref{eq1.10},
\begin{equation}
y(\lambda,x)=y_0s_0(-i\lambda x)+y_1\frac{s_1(-i\lambda x)}{(-i\lambda)}+y_2\frac{s_2(-i\lambda x)}{(-i\lambda)^2}+i\int\limits_0^x\frac{s_2(-i\lambda(x-t))}{(-i\lambda)^2}f(t)dt.\label{eq1.11}
\end{equation}
\vspace{5mm}

{\bf 1.2} Using differentiability of $p(x)$, re-write equation \eqref{eq1.1} as
\begin{equation}
y'''(\lambda,x)=i\lambda^3y(\lambda,x)+2p(x)y'(\lambda,x)+[p'(x)+iq(x)]y(\lambda,x)=0.\label{eq1.12}
\end{equation}
This equation, in view of conditions \eqref{eq1.3} for $|x|\rightarrow\infty$, becomes an elementary \eqref{eq1.4}, therefore, it is natural to define the Jost solutions $\{v_k(\lambda,x)\}_0^2$ and $\{u_k(\lambda,x)\}_0^2$ as solutions to equation \eqref{eq1.12} that asymptotically, for $|x|\rightarrow\infty$, behave as solutions to equation \eqref{eq1.4},
\begin{equation}
\begin{array}{ccc}
({\rm a})\quad v_k(\lambda,x)\rightarrow e^{-i\lambda\zeta_kx}&(x\rightarrow\infty,\,0\leq k\leq2);\\
({\rm b})\quad u_k(\lambda,x)\rightarrow e^{-i\lambda\zeta_kx}&(x\rightarrow-\infty,\,0\leq k\leq2).
\end{array}\label{eq1.13}
\end{equation}
Using \eqref{eq1.11}, for $v_k(\lambda,x)$, we obtain that
\begin{equation}
v_k(\lambda,x)=e^{-i\lambda\zeta_kx}-\int\limits_x^\infty\frac{s_2(-i\lambda(x-t))}{(-i\lambda)^2}\{2p(t)v'_k(\lambda,t)+[p'(t)+iq(t)]v_k(\lambda,t)\}dt\label{eq1.14}
\end{equation}
($0\leq k\leq2$) and thus
\begin{equation}
v'_k(\lambda,x)=-i\lambda\zeta_ke^{-i\lambda\zeta_kx}-\int\limits_x^\infty\frac{s_1(-i\lambda(x-t))}{(-i\lambda)}\{2p(t)v'_k(\lambda,t)+[p'(t)+iq(t)]v_k(\lambda,t)\}dt.\label{eq1.15}
\end{equation}
Multiplying \eqref{eq1.14} by $[p'(x)+iq(x)]$ and equality \eqref{eq1.15}, by $2p(x)$, and adding, we obtain for the function
\begin{equation}
w_k(\lambda,x)\stackrel{\rm def}{=}2p(x)v'_k(\lambda,x)+[p'(x)+iq(x)]v_k(\lambda,x),\label{eq1.16}
\end{equation}
equation
\begin{equation}
w_k(\lambda,x)=g_k(\lambda,x)-\int\limits_x^\infty K_1(\lambda,x,t)w_k(\lambda,t)dt\label{eq1.17}
\end{equation}
where
\begin{equation}
K_1(\lambda,x,t)=\frac{s_2(-i\lambda(x-t))}{(-i\lambda)^2}\left[p'(x)+q(x)+\frac{s_1(-i\lambda(x-t))}{-i\lambda}2p(x)\right];\label{eq1.18}
\end{equation}
\begin{equation}
g_k(\lambda,x)=e^{-i\lambda\zeta_kx}\{p'(x)+iq(x)-2i\lambda\zeta_kp(x)\}.\label{eq1.19}
\end{equation}

\begin{remark}\label{r1.2}
Equation \eqref{eq1.14} implies that $v_k(\lambda,x)$ is expressed via solution $w_k(\lambda,x)$ to equation \eqref{eq1.17} by the formula
\begin{equation}
v_k(\lambda,x)=e^{-i\lambda\zeta_kx}-\int\limits_x^\infty\frac{s_2(-i\lambda(x-t))}{(-i\lambda)^2}w_k(\lambda,t)dt.\label{eq1.20}
\end{equation}
\end{remark}

In the space $L^2(\mathbb{R})$, we define the family of Volterra operators
$$(K_\lambda f)(x)=\int\limits_x^\infty K_1(\lambda,x,t)f(t)dt\quad(f\in L^2(\mathbb{R})$$
where $K_1(\lambda,x,t)$ is given by \eqref{eq1.18}. In terms of $K_\lambda$, equation \eqref{eq1.17} becomes
$$(I+K_\lambda)w_k(\lambda,x)=g_k(\lambda,x),$$
and thus
\begin{equation}
w_k(\lambda,x)=\sum\limits_{n=0}^\infty(-1)^nK_\lambda^ng_k(\lambda,x).\label{eq1.21}
\end{equation}
Operators $K_\lambda^n$ are also Volterra,
$$(K_\lambda^nf)(x)=\int\limits_x^\infty K_n(\lambda,x,t)f(t)dt\quad(f\in L^2(\mathbb{R})),$$
and, for its kernels, the following recurrence relations hold:
\begin{equation}
K_{n+1}(\lambda,x,t)=\int\limits_x^tK_n(\lambda,x,s)K_n(\lambda,s,t)ds.\label{eq1.22}
\end{equation}
Since $|s_k(-i\lambda(x-t))|\leq e^{|\lambda|(t-x)}$ ($t\geq x$, $0\leq k\leq2$), then for $K_1(\lambda,x,t)$ \eqref{eq1.18}, the following inequality holds:
\begin{equation}
|K_1(\lambda,x,t)|\leq\frac1{|\lambda|^2}e^{|\lambda|(t-x)}m(\lambda,x)\label{eq1.23}
\end{equation}
where
\begin{equation}
m(\lambda,x)\stackrel{\rm def}{=}|p'(x)|+|q(x)|+2|\lambda||p(x)|.\label{eq1.24}
\end{equation}

\begin{lemma}\label{l1.1}
For kernels $K_n(\lambda,x,t)$, the following estimates hold:
\begin{equation}
|K_n(\lambda,x,t)|\leq\frac1{|\lambda|^{2n}}e^{|\lambda|(t-x)}\frac{m(\lambda,x)}{(n-1)!}\left(\int\limits_x^tm(\lambda,s)\right)^{n-1}\quad(t\geq x).\label{eq1.25}
\end{equation}
where $m(\lambda,x)$ is given by \eqref{eq1.24}.
\end{lemma}

P r o o f. Apply the method of mathematical induction by ``$n$''. For $n=1$, estimate \eqref{eq1.25} coincides with \eqref{eq1.23}. Using \eqref{eq1.22}, \eqref{eq1.23}, and \eqref{eq1.25} for ``$n$'', we obtain
$$|K_{n+1}(\lambda,x,t)|\leq\int\limits_x^t\frac1{|\lambda|^{2n}}e^{|\lambda|(s-x)}\frac{m(\lambda,x)}{(n-1)!}\left(\int\limits_x^sm(\lambda,\zeta)d\zeta\right)^{n-1}\frac1{|\lambda|^2}e^{|\lambda|(t-s)}m(\lambda,s)ds$$
$$=\frac1{|\lambda|^{2(n+1)}}e^{|\lambda|(t-x)}m(\lambda,x)\frac1{(n-1)!}\int\limits_x^t\left(\int\limits_x^sm(\lambda,\zeta)d\zeta\right)^{n-1}m(\lambda,s)ds$$
$$=\frac1{|\lambda|^{2(n+1)}}e^{|\lambda|(t-x)}m(\lambda,x)\frac1{n!}\left(\int\limits_x^tm(\lambda,s)ds\right)^n,$$
which coincides with \eqref{eq1.25} for ``$n+1$''. $\blacksquare$

Equation \eqref{eq1.21} implies that
\begin{equation}
w_k(\lambda,x)=g_k(\lambda,x)-\int\limits_x^\infty N(\lambda,x,t)g_k(\lambda,t)dt\label{eq1.26}
\end{equation}
where
\begin{equation}
N(\lambda,x,t)=\sum\limits_{n=1}^\infty(-1)^nK_n(\lambda,x,t).\label{eq1.27}
\end{equation}
Series \eqref{eq1.27}, in view of \eqref{eq1.25}, is majorized by the convergent series
$$|N(\lambda,x,t)|<\frac1{|\lambda|^2}e^{|\lambda|(t-x)}m(\lambda,x)\cdot\exp\left\{\frac1{|\lambda|}\cdot\int\limits_x^tm(\lambda,s)ds\right\}.$$
Note that
$$\frac1{|\lambda|^2}\int\limits_x^tm(\lambda,s)ds\leq\frac1{|\lambda|^2}\|p'(x)+iq(k)\|_{L^1(\mathbb{R})}+\frac2{|\lambda|}\|p\|_{L^1(\mathbb{R})}\stackrel{\rm def}{=}c(\lambda)$$
and $c(\lambda)\rightarrow0$ when $\lambda\rightarrow\infty$. So, for $N(\lambda,x,t)$ the following inequality holds:
\begin{equation}
|N(\lambda,x,t)|\leq\frac1{|\lambda|^2}e^{|\lambda|(t-x)}m(\lambda,x)e^{c(\lambda)}\quad(t>x,\lambda\not=0).\label{eq1.28}
\end{equation}
And since
$$|g_k(\lambda,x)|\leq e^{|\lambda|x}m(\lambda,x)\quad(x\geq0),$$
then, according to \eqref{eq1.26}, \eqref{eq1.28}, we have
$$|w_k(\lambda,x)|\leq e^{|\lambda|x}m(\lambda,x)+\frac1{|\lambda|^2}e^{c(\lambda)}\int\limits_x^\infty m(\lambda,x)e^{|\lambda|(t-x)}m(\lambda,t)e^{|\lambda|t}dt$$
$$=e^{|\lambda|x}m(\lambda,x)+\frac{m(\lambda,x)}{|\lambda|^2}e^{c(\lambda)}\int\limits_x^\infty e^{|\lambda|(2t-x)}m(\lambda,t)dt.$$
And since $e^{|\lambda|(2t-x)}m(\lambda,t)\leq e^{|\lambda|2t}m(\lambda,t)\in L^1(\mathbb{R})$, for $2|\lambda|\leq a$ (Remark \ref{r1.1}), then we obtain that
\begin{equation}
\begin{array}{ccc}
|w_k(\lambda,x)|\leq m(\lambda,x)\{e^{|\lambda|x}+M(\lambda)\}\quad(\forall\lambda\not=0,2|\lambda|<a);\\
M(\lambda)\stackrel{\rm def}{=}e^{c(\lambda)}\|e^{2|\lambda|t}m(\lambda,t)\|_{L^1(\mathbb{R})}.
\end{array}\label{eq1.29}
\end{equation}
Hence it follows that $v_k(\lambda,x)$ as a solution to the integral equation \eqref{eq1.14} satisfies the boundary condition (a) \eqref{eq1.13}. Really, using \eqref{eq1.20}, we obtain
$$\left|v_k(\lambda,x)-e^{-i\lambda\zeta_kx}\right|\leq\int\limits_x^\infty\left|\frac{s_2(-i\lambda(x-t))}{(-i\lambda)^2}\right||w_k(\lambda,t)|dt$$
and, taking into account \eqref{eq1.29}, we have
\begin{equation}
\left|v_k(\lambda,x)-e^{-i\lambda\zeta_kx}\right|\leq\frac1{|\lambda|^2}\int\limits_x^\infty e^{|\lambda|(t-x)}m(\lambda,t)\left\{e^{|\lambda|t}+M(\lambda)\right\}dt.\label{eq1.30}
\end{equation}
And since $e^{|\lambda|(2t-x)}m(\lambda,t)\leq e^{2|\lambda|t}m(\lambda,t)\in L^1(\mathbb{R})$, for $2|\lambda|\leq a$, and $e^{|\lambda|(t-x)}m(\lambda,t)\leq e^{|\lambda|t}m(\lambda,t)\in L^1(\mathbb{R})$, for $|\lambda|\leq a$ (Remark \ref{r1.1}), then integral in \eqref{eq1.30} vanishes when $x\rightarrow\infty$. Define the disk of the radius $a/2$,
\begin{equation}
\mathbb{D}_0(a/2)=\left\{\lambda\in\mathbb{C}:\lambda\not=0,|\lambda|<\frac a2\right\}.\label{eq1.31}
\end{equation}

\begin{theorem}\label{t1.1}
If conditions \eqref{eq1.3} hold, then for all $\lambda\in\mathbb{D}_0(a/2)$ \eqref{eq1.31}, Jost solutions $\{v_k(\lambda,x)\}_0^2$ to the boundary problem \eqref{eq1.12}, {\rm(a)} \eqref{eq1.13} exist and are given by \eqref{eq1.20} where $w_k(\lambda,x)$ are solutions to the integral equation \eqref{eq1.17}. Functions $\{v_k(\lambda,x)\}_0^2$ and $\{w_k(\lambda,x)\}_0^2$ are holomorphic inside the disk $\mathbb{D}_0(a/2)$ \eqref{eq1.31}.
\end{theorem}

\begin{remark}\label{r1.3}
For the functions $v_k(\lambda,x)$ and $w_k(\lambda,x)$, the following equalities hold:
\begin{equation}
v_k(\lambda\zeta_1,x)=v_{k'}(\lambda,x);\quad w_k(\lambda\zeta_1,x)=w_{k'}(\lambda,x);\quad k'=(k+1)(\hspace{-6mm}\mod3).\label{eq1.32}
\end{equation}
\end{remark}
\vspace{5mm}

{\bf 1.3} Analogously to \eqref{eq1.14}, integral equation
\begin{equation}
u_k(\lambda,x)=e^{-i\lambda\zeta_kx}+\int\limits_{-\infty}^x\frac{s_2(-i\lambda(x-t))}{(-i\lambda)^2}\{2p(t)u'_k(\lambda,t)+[p(t)+iq(t)]u_k(\lambda,t)\}dt\label{eq1.33}
\end{equation}
($0\leq k\leq2$) corresponds to the Jost solutions $\{u_k(\lambda,x)\}$ to the boundary value problem \eqref{eq1.12}, (b) \eqref{eq1.13}, and thus
$$u'_k(\lambda,x)=-i\lambda\zeta_ke^{-i\lambda\zeta_kx}+\int\limits_{-\infty}^x\frac{s_1(-i\lambda(x-t))}{(-i\lambda)}\{2p(t)u'_k(\lambda,t)+[p(t)+iq(t)]u_k(\lambda,t)\}dt.$$
Hence, for the functions
\begin{equation}
\widehat{w}_k(\lambda,x)=2p(x)u'_k(\lambda,x)+[p(x)+iq(x)]u_k(\lambda,x)\quad(0\leq k\leq2),\label{eq1.34}
\end{equation}
we obtain the integral Volterra equation
\begin{equation}
\widehat{w}_k(\lambda,x)=g_k(\lambda,x)+\int\limits_{-\infty}^xK_1(\lambda,x,t)\widehat{w}_k(\lambda,t)dt\quad(0\leq k\leq2)\label{eq1.35}
\end{equation}
where $K_1(\lambda,x,t)$ and $g_k(\lambda,x)$ are given by \eqref{eq1.18} and \eqref{eq1.19} correspondingly.

\begin{remark}\label{r1.4}
Knowing the solution $\{\widehat{w}_k(\lambda,x)\}_0^2$ to the integral equation \eqref{eq1.35}, in view of \eqref{eq1.33}, we obtain solutions $\{u_k(\lambda,x)\}_0^2$ to equation \eqref{eq1.33},
\begin{equation}
u_k(\lambda,x)=e^{-i\lambda\zeta_kx}+\int\limits_{-\infty}^x\frac{s_2(-i\lambda(x-t))}{(-i\lambda)^2}\widehat{w}_k(\lambda,t)dt.\label{eq1.36}
\end{equation}
\end{remark}

For Jost solutions $\{u_k(\lambda,x)\}_0^2$, an analogue to Theorem \ref{t1.1} is true.

\begin{theorem}\label{t1.2}
If conditions \eqref{eq1.3} are met, then for all $\lambda\in\mathbb{D}_0(a/2)$ \eqref{eq1.31}, Jost solutions $\{u_k(\lambda,x)\}_0^2$ exist and are given by \eqref{eq1.36} where $\{\widehat{w}_k(\lambda,x)\}_0^2$ are solutions to the integral equation \eqref{eq1.35}. Functions $\{u_k(\lambda,x)\}_0^2$ and $\{\widehat{w}_k(\lambda,x)\}_0^2$ are analytical in the disk $\mathbb{D}_0(a/2)$ \eqref{eq1.31}.
\end{theorem}

For functions $\{u_k(\lambda,x)\}_0^2$ and $\{\widehat{w}_k(\lambda,x)\}_0^2$, equalities \eqref{eq1.32} hold.
\vspace{5mm}

{\bf 1.4} Proceed to analytical (with respect to $\lambda$) properties of the functions $\{v_k(\lambda,x)\}_0^2$ and $\{u_k(\lambda,x)\}_0^2$. Define the functions
\begin{equation}
\psi_k(\lambda,x)\stackrel{\rm def}{=}v_k(\lambda,x)e^{i\lambda\zeta_kx}\quad(0\leq k\leq2),\label{eq1.37}
\end{equation}
then \eqref{eq1.20} implies
\begin{equation}
\psi_k(\lambda,x)=1-\int\limits_x^\infty e^{i\lambda(x-t)\zeta_k}\frac{s_2(-i\lambda(x-t))}{(-i\lambda)^2}{\rm v}_k(\lambda,t)dt,\label{eq1.38}
\end{equation}
besides, functions
\begin{equation}
{\rm v}_k(\lambda,x)\stackrel{\rm def}{=}w_k(\lambda,x)e^{i\lambda\zeta_kx},\label{eq1.39}
\end{equation}
in its turn, due to \eqref{eq1.17}, are solutions to the equations
\begin{equation}
{\rm v}_k(\lambda,x)=f_k(\lambda,x)-\int\limits_x^\infty e^{i\lambda\zeta_k(x-t)}K_1(\lambda,x,t){\rm v}_k(\lambda,t)dt\label{eq1.40}
\end{equation}
where
\begin{equation}
f_k(\lambda,x)\stackrel{\rm def}{=}g_k(\lambda,x)e^{i\lambda\zeta_kx}=p'(x)+iq(x)-2i\lambda\zeta_kp(x).\label{eq1.41}
\end{equation}
Consider the case of $k=0$. Since
\begin{equation}
e^{i\lambda\zeta_0(x-t)}s_p(-i\lambda(x-t))=\frac13\left\{1+\frac1{\zeta_1^p}e^{i\lambda(\zeta_0-\zeta_1)(t-x)}+\frac1{\zeta_2^p}e^{i\lambda(\zeta_2-\zeta_0)(t-x)}\right\}\label{eq1.42}
\end{equation}
($0\leq p\leq2$), and
$$i\lambda(\zeta_1-\zeta_0)=\frac12[(3\beta-\alpha\sqrt3)+i(-\beta\sqrt3-\alpha3)];\quad i\lambda(\zeta_2-\zeta_0)=\frac12[(3\beta+\alpha\sqrt3)+i(\beta\sqrt3-3\alpha)]$$
($\lambda=\alpha+i\beta\in\mathbb{C}$, $\alpha$, $\beta\in\mathbb{R}$), then
$$\left|e^{i\lambda(\zeta_1-\zeta_0)(t-x)}\right|=e^{\frac12(3\beta-\sqrt3\alpha)(t-x)};\quad\left|e^{i\lambda(\zeta_2-\zeta_0)(t-x)}\right|=e^{\frac12(3\beta+\alpha\sqrt3)(t-x)}.$$
And thus, for $\beta\sqrt3-\alpha<0$, $\beta\sqrt3<-\alpha$, when $t\geq x$, each of these exponents is lesser than one, and for $\lambda\rightarrow$ in this domain these exponents decrease in modulo and vanish. Domain $\{\lambda\in\mathbb{C}:\beta\sqrt3<\alpha,\beta\sqrt3<-\alpha\}$ is convenient to rewrite in terms of sectors $\{S_p\}_0^5$ \eqref{eq1.9}. By $\{S_p(i)\}_0^5$, we denote the sectors obtained from $\{S_p\}_0^5$ after the rotation through $\pi/2$,
\begin{equation}
S_p(i)\stackrel{\rm def}{=}\{i\lambda:\lambda\in S_p\}\quad(0\leq p\leq5).\label{eq1.43}
\end{equation}

\begin{picture}(200,200)
\put(0,100){\line(1,0){200}}
\put(100,0){\line(0,1){200}}
\put(0,133){\line(3,-1){200}}
\put(200,133){\line(-3,-1){200}}
%\put(105,185){$il_{\zeta_0}$}
\put(120,150){$\Omega_2$}
\qbezier(153,120)(100,140)(49,120)
\qbezier(100,23)(155,35)(168,125)
\qbezier(45,120)(30,70)(100,25)
%\qbezier(125,100)(115,110)(120,104)
\put(140,105){$S_4(i)$}
\put(70,110){$S_0(i)$}
%\put(80,88){$S_4$}
\put(70,70){$S_2(i)$}
\put(50,100){$S_1(i)$}
\put(30,140){$\Omega_1$}
%\put(0,50){$il_{\zeta_1}$}
\put(115,70){$S_3(i)$}
\put(105,110){$S_5(i)$}
\put(180,60){$\Omega_0$}
%\put(185,105){$l_{\zeta_0}$}
%\put(180,60){$il_{\zeta_2}$}
%\put(165,50){$\Omega_1$}
%\put(80,40){$\psi_1^*(\lambda,0)$}
\end{picture}

\hspace{20mm} Fig. 1

Using $\{S_p(i)\}_0^5$ \eqref{eq1.43}, we define three sectors $\{\Omega_p\}_0^2$ in $\mathbb{C}$:
$$\Omega_0\stackrel{\rm def}{=}S_2(i)\cup S_3(i)\cup(i\widehat{l_{\zeta_0}});$$
\begin{equation}
\Omega_2\stackrel{\rm def}{=}S_4(i)\cup S_5(i)\cup(i\widehat{l}_{\zeta_1});\label{eq1.44}
\end{equation}
$$\Omega_1\stackrel{\rm def}{=}S_0(i)\cup S_1(i)\cup(i\widehat{l}_{\zeta_2}).$$
It is easy to see that the domain $\{\lambda\in\mathbb{C}:\sqrt3\beta<\alpha,\sqrt3\beta<-\alpha\}$ coincides with the sector $\Omega_0$. These considerations imply that
\begin{equation}
e^{i\lambda\zeta_0(x-t)}s_p(-i\lambda(x-t))\rightarrow\frac13\quad(\lambda\rightarrow\infty,\lambda\in\Omega_0,x>t,0\leq p\leq2).\label{eq1.45}
\end{equation}
Using \eqref{eq1.18}, we have that
\begin{equation}
e^{i\lambda\zeta_0(x-t)}K_1(\lambda,x,t)\sim\frac1{3(i\lambda)^2}(p'(x)+iq(x))+\frac2{3(i\lambda)}p(x)\quad(\lambda\in\Omega_0,\,|\lambda|\gg1).\label{eq1.46}
\end{equation}
Therefore, function ${\rm v}_0(\lambda,x)$, due to equation \eqref{eq1.40}, for $|\lambda|\gg1$ and $\lambda\in\Omega_0$ behaves in the following way:
\begin{equation}
{\rm v}_0(\lambda,x)\sim p'(x)+iq(x)+2ip(x)+o\left(\frac1{\lambda^2}\right)\quad(|\lambda|\gg1,\lambda\in\Omega_0).\label{eq1.47}
\end{equation}
As a result, we arrive at the statement
\begin{theorem}\label{t1.3}
Functions $v_0(\lambda,x)$ and $\psi_0(\lambda,x)$ \eqref{eq1.37} are holomorphic in the sector $\Omega_0$ and the following relations hold:
\begin{equation}
\begin{array}{lll}
{\displaystyle({\rm i})\,\lim\limits_{\stackrel{\lambda\rightarrow\infty}{\lambda\in\Omega_0}}(\psi_0(\lambda,x)-1)(-i\lambda)^2\frac3{2i\lambda}=P(x);\quad P(x)\stackrel{\rm def}{=}\int\limits_x^\infty p(t)dt;}\\
{\displaystyle({\rm ii})\,\lim\limits_{\stackrel{\lambda\rightarrow\infty}{\lambda\in\Omega_0}}[\psi_0(\lambda,x)-1+2i\lambda P(x)]3(-i\lambda)^2=Q(x);\quad Q(x)\stackrel{\rm def}{=}-\int\limits_x^\infty[p'(t)+iq(t)]dt.}
\end{array}\label{eq1.48}
\end{equation}
\end{theorem}

So, by the function $\psi_0(\lambda,x)$ we find  functions $P(x)$ and $Q(x)$ whence, upon differentiating, we recover potentials on the right half-axis,
\begin{equation}
P'(x)=p(x)\quad Q'(x)=(p'(x)+iq(x))\quad(x\in\mathbb{R}_+).\label{eq1.49}
\end{equation}
Equation \eqref{eq1.32} implies that functions $\psi_1(\lambda,x)$, $\psi_2(\lambda,x)$ \eqref{eq1.37} are matched with sectors $\Omega_1$ and $\Omega_2$ \eqref{eq1.44} correspondingly, and analogues of Theorem \ref{t1.3} are true for these functions in the given sectors.
\vspace{5mm}

{\bf 1.5} Analogously to \eqref{eq1.37}, Jost solutions $\{u_k(\lambda,x)\}_0^2$ with
\begin{equation}
\varphi_k(\lambda,x)=u_k(\lambda,x)e^{i\lambda\zeta_kx}\quad(0\leq k\leq2),\label{eq1.50}
\end{equation}
then \eqref{eq1.35} implies that
\begin{equation}
\varphi_k(\lambda,x)=1+\int\limits_{-\infty}^xe^{i\lambda\zeta_k(x-t)}\frac{s_2(-i\lambda(x-t))}{(-i\lambda)^2}\widehat{\rm v}_k(\lambda,t)dt,\label{eq1.51}
\end{equation}
besides, functions
\begin{equation}
\widehat{\rm v}_k(\lambda,x)=\widehat{\rm w}_k(\lambda,x)e^{i\lambda\zeta_kx},\label{eq1.52}
\end{equation}
due to \eqref{eq1.34}, are solutions to the equations
\begin{equation}
\widehat{\rm v}_k(\lambda,x)=f_k(\lambda,x)+\int\limits_{-\infty}^xe^{i\lambda\zeta_k(x-t)}K_1(\lambda,x,t)\widehat{\rm v}_k(\lambda,t)dt\label{eq1.53}
\end{equation}
where $K_1(\lambda,x,t)$ is given by \eqref{eq1.18} and $f_k(\lambda,x)$, by \eqref{eq1.41}.

By $\{\Omega_k^-\}_0^2$, we denote the centrally symmetrical sectors $\{\Omega_k\}_0^2$ \eqref{eq1.44},
\begin{equation}
\Omega_k^-=\{\lambda\in\mathbb{C}:-\lambda\in\Omega_k\}\quad(0\leq k\leq2).\label{eq1.54}
\end{equation}
Analogously to Theorem \ref{t1.3}, the following statement is true.

\begin{theorem}\label{t1.4}
Functions $u_0(\lambda,x)$ and $\varphi_0(\lambda,x)$ \eqref{eq1.50} are holomorphic in the sector $\Omega_0^-$ and
\begin{equation}
\begin{array}{lll}
{\displaystyle({\rm i})\quad\lim\limits_{\lambda\rightarrow\infty\atop\lambda\in\Omega_0^-}(\varphi_0(\lambda,x)-1)(-i\lambda)^2\frac3{2i\lambda}=\widehat{P}(x);\quad\widehat{P}(x)\stackrel{\rm def}{=}\int\limits_{-\infty}^xp(t)dt;}\\
{\displaystyle({\rm ii})\quad\lim\limits_{\lambda\rightarrow\infty\atop\lambda\in\Omega_0^-}(\varphi_0(\lambda,x)-1-2i\lambda \widehat{P}(x))3(-i\lambda)^2=\widehat{Q}(x);\quad\widehat{Q}(x)=\int\limits_{-\infty}^x[p'(t)+iq(t)]dt.}
\end{array}\label{eq1.55}
\end{equation}
\end{theorem}

Thus, from the function $\varphi_0(\lambda,x)$, we find $\widehat{P}(x)$ and $\widehat{Q}(x)$, upon differentiating it allows us to find potentials on the left half-axis,
\begin{equation}
\widehat{P}'(x)=p(x);\quad\widehat{Q}'=p'(x)+iq(x)\quad(x\in\mathbb{R}_-)\label{eq1.56}
\end{equation}
\vspace{5mm}

{\bf 1.6} Study the problem of linear independence of Jost solutions $\{v_k(\lambda,x)\}_0^2$ (and $\{u_k(\lambda,x)\}_0^2$). Consider the Wronskian
\begin{equation}
\Delta(\lambda,x)=\det\left[
\begin{array}{ccc}
v_0(\lambda,x)&v_1(\lambda,x)&v_2(\lambda,x)\\
v'_0(\lambda,x)&v'_1(\lambda,x)&v'_2(\lambda,x)\\
v''_0(\lambda,x)&v''_1(\lambda,x)&v''_2(\lambda,x)
\end{array}\right].\label{eq1.57}
\end{equation}
Evidently, $\Delta'(\lambda,x)=0$, due to equation \eqref{eq1.12}. Using (a) \eqref{eq1.13}, we obtain
$$\Delta(\lambda,x)=\det\left[
\begin{array}{ccc}
e^{-i\lambda\zeta_0x}&e^{-i\lambda\zeta_1x}&e^{-i\lambda\zeta_2x}\\
(-i\lambda\zeta_0)e^{-i\lambda\zeta_0x}&(-i\lambda\zeta_1)e^{-i\lambda\zeta_1x}&(-i\lambda\zeta_2)e^{-i\lambda\zeta_2x}\\
(-i\lambda\zeta_0)^2e^{-i\lambda\zeta_0x}&(-i\lambda\zeta_1)^2e^{-i\lambda\zeta_0x}&(-i\lambda\zeta_1)^2e^{-i\lambda\zeta_1x}
\end{array}\right],$$
and since $\zeta_0+\zeta_1+\zeta_2=0$, then
$$\Delta(\lambda,x)=(-i\lambda)^3\det\left[
\begin{array}{ccc}
1&1&1\\
\zeta_0&\zeta_1&\zeta_2\\
\zeta_0^2&\zeta_1^2&\zeta_1^2
\end{array}\right]=(-i\lambda)^3\det\left[
\begin{array}{ccc}
1&1&1\\
1&\zeta_1&\zeta_2\\
1&\zeta_2&\zeta_1
\end{array}\right]$$
\begin{equation}
=(-i\lambda)^3\cdot i\sqrt3=\sqrt3\lambda^3.\label{eq1.57}
\end{equation}

\begin{lemma}
For all $\lambda\not=0$, Jost solutions $\{v_k(\lambda,x)\}_0^2$ (correspondingly, $\{u_k(\lambda,x)\}_0^2$) are linearly independent.
\end{lemma}

\section{Scattering problem (waves incident from $+\infty$)}

{\bf 2.1} Decompose every Jost solution $\{u_k(\lambda,x)\}_0^2$ into other Jost solutions $\{v_k(\lambda,x)\}_0^2$,
\begin{equation}
u_k(\lambda,x)=\sum\limits_lt_{k,l}(\lambda)v_l(\lambda,x)\quad(0\leq k\leq2),\label{eq2.1}
\end{equation}
or
\begin{equation}
u(\lambda,x)=T(\lambda)v(\lambda,x)\label{eq2.2}
\end{equation}
where $u(\lambda,x)\stackrel{\rm def}{=}\col[u_0(\lambda,x),u_1(\lambda,x),u_2(\lambda,x)]$; $v(\lambda,x)\stackrel{\rm def}{=}\col[v_0(\lambda,x),v_1(\lambda,x),$ $v_2(\lambda,x)]$, and $T(\lambda)$ is the {\bf transition matrix},
\begin{equation}
T(\lambda)\stackrel{\rm def}{=}\left[
\begin{array}{ccc}
t_{0,0}(\lambda)&t_{0,1}(\lambda)&t_{0,2}(\lambda)\\
t_{1,0}(\lambda)&t_{1,1}(\lambda)&t_{1,2}(\lambda)\\
t_{2,0}(\lambda)&t_{2,1}(\lambda)&t_{2,2}(\lambda)
\end{array}\right].\label{eq2.3}
\end{equation}

\begin{remark}\label{r2.1}
It is sufficient to consider only one equality among \eqref{eq2.1}, e.g., for $k=0$,
\begin{equation}
u_0(\lambda,x)=\sum\limits_lt_{0,l}(\lambda)v_l(\lambda,x)\label{eq2.4}
\end{equation}
all the other follow from it since $u_0(\lambda\zeta_1)=u_1(\lambda,x)$, $u_0(\lambda\zeta_2)=u_2(\lambda,x)$, and $v_0(\lambda\zeta_1)=v_1(\lambda,x)$, $v_0(\lambda\zeta_2,x)=v_2(\lambda,x)$. Besides,
\begin{equation}
\begin{array}{ccc}
t_{0,0}(\lambda\zeta_1)=t_{0,1}(\lambda);&t_{0,1}(\lambda\zeta_1)=t_{1,1}(\lambda);&t_{0,2}(\lambda\zeta_1)=t_{1,0}(\lambda);\\
t_{0,0}(\lambda\zeta_2)=t_{2,2}(\lambda);&t_{0,1}(\lambda\zeta_2)=t_{2,0}(\lambda);&t_{0,2}(\lambda\zeta_2)=t_{2,1}(\lambda).
\end{array}\label{eq2.5}
\end{equation}
\end{remark}

By $W_{k,s}(v,\lambda,x)$, we denote the Wronskian of functions $v_k(\lambda,x)$ and $v_s(\lambda,x)$,
\begin{equation}
W_{k,s}(v,\lambda,x)=v'_k(\lambda,x)v_s(\lambda,x)-v'_s(\lambda,x)v_k(\lambda,x)\quad(0\leq k,s\leq2).\label{eq2.6}
\end{equation}
Define the involution ``$+$'',
\begin{equation}
f^+(\lambda)\stackrel{\rm def}{=}\overline{f(\overline{\lambda})}.\label{eq2.7}
\end{equation}

\begin{lemma}\label{l2.1}
For Wronskians $W_{k,s}(v,\lambda,x)$, the following formulas are true:
\begin{equation}
\begin{array}{ccc}
W_{1,2}(v,\lambda,x)=\sqrt3\lambda\zeta_0v_0^+(\lambda,x);\,W_{2,0}(v,\lambda,x)=\sqrt3\lambda\zeta_1v_1^+(\lambda,x);\\
W_{0,1}(v,\lambda,x)=\sqrt3\lambda\zeta_2v_2^+(\lambda,x).
\end{array}\label{eq2.8}
\end{equation}
\end{lemma}

P r o o f. For $W_{k,s}(v)=v'_kv_s-v'_sv_k$, we have
$$W'_{k,s}(v)=v''_kv_s-v''_kv_k,$$
whence it follows that
$$W''_{k,s}(v)=v'''_kv_s-v'''_sv_k+v''_kv'_s-v''_sv'_k.$$
Use the fact that $\{v_k\}_0^2$ are solutions to equation
\begin{equation}
y'''=i\lambda^3y+2py'=(p+iq)y,\label{eq2.9}
\end{equation}
then we have
$$W''_{k,s}=(i\lambda^3v_k+2pv'_k+(p+iq)v_k)v_s+(i\lambda^3v_s+2pv'_s+(p+iq)v_s))v_k$$
$$+v''_kv'_s-v''_sv'_k=2pW_{k,s}+iv''_kv'_s-v''_sv'_k.$$
Upon differentiating this equality again, we obtain
$$W'''_{k,s}=2p'W_{k,s}+2pW'_{k,s}+(i\lambda^3v_k+2pv'_k+(p+iq)v_k)v'_s$$
$$-(i\lambda^3v_s+2pv'_s+(p+iq)v_s)v'_k=-i\lambda^3W_{k,s}+2pW'_{k,s}+(p-iq)W_{k,s},$$
and thus $W_{k,s}$ is the solution to equation
\begin{equation}
y'''=-i\lambda^3y+2py'+(p-ia)y\label{eq2.10}
\end{equation}
which we can derive from equation \eqref{eq2.9} upon substituting $\lambda\rightarrow\overline{\lambda}$ followed by conjugation of coefficients. In order to obtain formulas \eqref{eq2.8}, we need to calculate asymptotic of $W_{k,s}(v,\lambda,x)$ as $x\rightarrow\infty$ since for $W_{1,2}(v,\lambda,x)\rightarrow i\lambda(\zeta_2-\zeta_1)$ $e^{-i\lambda(\zeta_1+\zeta_2)}$. It is left to note that $\zeta_2-\zeta_1=-i\sqrt3$ and $\zeta_2+\zeta_3+\zeta_0=0$. $\blacksquare$

\begin{remark}\label{r2.2}
For Wronskians
\begin{equation}
W_{k,s}(u,\lambda,x)\stackrel{\rm def}{=}u'_k(\lambda,x)u_s(\lambda,x)-u'_s(\lambda,x)u_k(\lambda,x)\label{eq2.11}
\end{equation}
analogues of the formulas \eqref{eq2.8} hold.
\end{remark}
\vspace{5mm}

{\bf 2.2} Using \eqref{eq2.1}, calculate Wronskian of the functions $u_1(\lambda,x)$ and $u_2(\lambda,x)$,
$$W_{1,2}(u,\lambda,x)=W\{(t_{1,0}(\lambda)v_0(\lambda)+t_{1,1}(\lambda)v_1(\lambda)+t_{1,2}(\lambda)v_2(\lambda)),(t_{2,0}(\lambda)v_0(\lambda,x)$$
$$+t_{2,1}(\lambda)v_1(\lambda,x)+t_{2,2}(\lambda)v_2(\lambda,x))\}=W_{0,1}(v,\lambda,x)\left|
\begin{array}{ccc}
t_{1,0}(\lambda)&t_{1,1}(\lambda)\\
t_{2,0}(\lambda)&t_{2,1}(\lambda)
\end{array}\right|$$
$$+W_{0,2}(v,\lambda,x)\left|
\begin{array}{ccc}
t_{1,0}(\lambda)&t_{1,2}(\lambda)\\
t_{2,0}(\lambda)&t_{2,2}(\lambda)
\end{array}\right|+W_{1,2}(v,\lambda,x)\left|
\begin{array}{ccc}
t_{1,1}(\lambda)&t_{1,2}(\lambda)\\
t_{2,1}(\lambda)&t_{2,2}(\lambda)
\end{array}\right|.$$
In view of \eqref{eq2.8}, we obtain
$$u_0^+(\lambda,x)=v_0^+(\lambda,x)T_{0,0}(\lambda)+\zeta_1v_1^+(\lambda,x)T_{0,1}(\lambda)+\zeta_2v_2^+(\lambda,x)T_{0,2}(\lambda)$$
where $T_{k,s}(\lambda)$ are algebraic complements of the elements $t_{k,s}(\lambda)$ of the matrix $T(\lambda)$ \eqref{eq2.3}. Analogously, we find that
$$\zeta_2u_2^+(\lambda,x)=v_0^+(\lambda,x)T_{2,0}(\lambda)+\zeta_1v_1^+(\lambda,x)T_{2,1}(\lambda)+\zeta_2v_2^+(\lambda,x)T_{2,2}(\lambda);$$
$$\zeta_1u_1^+(\lambda,x)=v_0^+(\lambda,x)T_{1,0}(\lambda)+\zeta_1v_1^+(\lambda,x)T_{1,1}(\lambda)+\zeta_2v_2^+(\lambda,x)T_{1,2}(\lambda).$$
Define the matrices
\begin{equation}
\widehat{T}(\lambda)\stackrel{\rm def}{=}\left[
\begin{array}{ccc}
T_{0,0}(\lambda)&T_{0,1}(\lambda)&T_{0,2}(\lambda)\\
T_{1,0}(\lambda)&T_{1,1}(\lambda)&T_{1,2}(\lambda)\\
T_{2,0}(\lambda)&T_{2,1}(\lambda)&T_{2,2}(\lambda)
\end{array}\right];\quad J=\left[
\begin{array}{ccc}
1&0&0\\
0&\zeta_1&0\\
0&0&\zeta_2
\end{array}\right],\label{eq2.12}
\end{equation}
then the obtained equalities in the matrix form are
\begin{equation}
T(\lambda)Jv^+(\lambda,x)=Ju^+(\lambda,x)\label{eq2.13}
\end{equation}
where $u^+(\lambda,x)$ and $v^+(\lambda,x)$ are obtained from $u(\lambda,x)$ and $v(\lambda,x)$ upon applying the involution ``$+$'' \eqref{eq2.7}. Since $\widehat{T}(\lambda)=\Delta(\lambda)\left(T^t(\lambda)\right)^{-1}$ where $\Delta(\lambda)=\det T(\lambda)$ and $T^t(\lambda)$ is a matrix obtained from $T(\lambda)$ by transposition. Therefore, due to \eqref{eq2.13}, we obtain
$$\Delta(\lambda)Jv^+(\lambda)=T^t(\lambda)Ju^+(\lambda,x).$$
Upon applying involution ``+'' \eqref{eq2.7} to both sides of equality, we obtain
$$\Delta^+(\lambda)J^*v(\lambda)=T^*(\overline{\lambda})J^*u(\lambda,x).$$
Upon applying $J$ ($JJ^*=I$) to both sides, we have
$$\Delta^+(\lambda)v(\lambda,x)=JT^*(\overline{\lambda})J^*u(\lambda,x),$$
and, due to \eqref{eq2.2},
$$\Delta^+(\lambda)u(\lambda,x)=T(\lambda)JT^*(\overline{\lambda})J^*u(\lambda,x),$$
which gives us
\begin{equation}
\Delta^+(\lambda)J=T(\lambda)JT^*(\overline{\lambda}).\label{eq2.14}
\end{equation}

\begin{lemma}\label{l2.2}
Matrix $T(\lambda)$ have the properties
\begin{equation}
\begin{array}{lll}
({\rm a})\quad\det T(\lambda)=\overline{\zeta}_p;\\
({\rm b})\quad\zeta_pJ=T(\lambda)JT^*(\overline{\lambda}),
\end{array}\label{eq2.15}
\end{equation}
where $\zeta_p$ is one of the roots $\zeta_k\}_0^2$ \eqref{eq1.5}.
\end{lemma}

P r o o f. Calculate determinants in both sides of equality \eqref{eq2.14}. Then we obtain $(\Delta^+(\lambda))^3=\Delta(\lambda)\Delta^+(\lambda)$ ($\det J=1$). Apply involution ``+'' \eqref{eq1.7} to both sides of this equality, $\Delta^3(\lambda)=\Delta^+(\lambda)\Delta(\lambda)$ and subtract it from the initial one, as a result we have $(\Delta^+(\lambda))^3=\Delta(\lambda)$, i.e., $(\Delta(\lambda)/\Delta^+(\lambda))^3=1$, and thus $\Delta(\lambda)=\zeta_q\Delta^+(\lambda)$, $\zeta_q$ is one of the roots $\{\zeta_k\}_0^2$ \eqref{eq1.5}. Upon substituting this expression into the initial equality, we obtain $(\Delta^+(\lambda))^2(\Delta^+(\lambda)-\zeta_q)=0$. So, either $\Delta^+(\lambda)=0$ or $\Delta^+(\lambda)=\zeta_q$. Equality $\Delta^+(\lambda)=0$ contradicts the linear independence of $\{u_k(\lambda,x)\}_0^2$. Therefore, $\Delta^+(\lambda)=\zeta_q$, which concludes the proof ($\zeta_p=\overline{\zeta_q}$). $\blacksquare$
\vspace{5mm}

{\bf 2.3} Proceed to bound states of the operator $L_{p,q}$. Its self-adjointness implies that $\lambda^3$ in equation \eqref{eq1.12} is real and thus points $\lambda$ lie on the straight lines $\{L_{\zeta_k}\}_0^2$ \eqref{eq1.7}. If $\mu>0$ from $\Omega_0^-$ is a common zero of the functions $t_{0,0}(\lambda)$ and $t_{0,1}(\lambda)$, then
$$u_0(\mu,x)=t_{0,2}(\mu)v_2(\mu,x).$$
Since $u_0(\mu,x)\in L^2(\mathbb{R})$ and $v_2(\mu,x)\in L^2(\mathbb{R}_+)$, then $u_0(\mu,x)\in L^2(\mathbb{R})$ is a bound state of the operator $L_{p,q}$.

Analogously, if $\nu<0$ from $\Omega_0^-$ is a common zero of $t_{0,0}(\lambda)$ and $t_{0,2}(\lambda)$, then $u_0(\nu,x)=t_{0,1}(\nu)v_1(\nu,x)$ is also a bound state of the operator $L_{p,q}$. Poits $\mu\zeta_1$, $\nu\zeta_2$ belong to the rays $\mu\zeta_1\in l_{\zeta_1}$ and $\nu\zeta_2\in\widehat{l}_{\zeta_2}$. Analogous position of zeros takes place in the sectors $\Omega_1^-$ and $\Omega_2^-$, which corresponds to the functions $t_{1,1}(\lambda)$ and $t_{2,2}(\lambda)$ correspondingly.

To calculate $t_{0,0}(\lambda)$, consider the system of equations that follows from \eqref{eq2.4},
$$\left\{
\begin{array}{lll}
u_0(\lambda,x)=t_{0,0}(\lambda)v_0(\lambda,x)+t_{0,1}(\lambda)v_1(\lambda,x)+t_{0,2}(\lambda)v_2(\lambda,x);\\
u'_0(\lambda,x)=t_{0,0}(\lambda)v'_0(\lambda,x)+t_{0,1}(\lambda)v'_1(\lambda,x)+t_{0,2}(\lambda)v'_2(\lambda,x);\\
u''_0(\lambda,x)=t_{0,0}(\lambda)v''_0(\lambda,x)+t_{0,1}(\lambda)v''_1(\lambda,x)+t_{0,2}(\lambda)v''_2(\lambda,x).
\end{array}\right.$$
Determinant of this system $\Delta(\lambda,x)=\sqrt3\lambda^3$ \eqref{eq1.57}. From this system of equations, we find
\begin{equation}
\begin{array}{ccc}
{\displaystyle t_{0,0}(\lambda)=\frac1{\sqrt3\lambda^3}\{u_0(\lambda,x)W''_{2,1}(\lambda,x)-u'_0(\lambda,x)W_{2,1}(\lambda,x)+u''_0(\lambda,x)W_{2,1}(\lambda,x)\}}\\
{\displaystyle=-\frac1{\lambda^2}\{u_0(\lambda,x)(e_0^+(\lambda,x))''-u'_0(\lambda,x)(e_0^+(\lambda,x))'+u''_0(\lambda,x)e_0^+(\lambda,x)\}.}
\end{array}\label{eq2.16}
\end{equation}

It is easy to check that the right-hand side of equality \eqref{eq2.16} does not depend on $x$. Therefore, we assume that $x=0$.

\begin{lemma}\label{l2.3}
Function $t_{0,0}(\lambda)$ \eqref{eq2.16} is holomorphic in the sector $\Omega_0^-$ and has there a finite number of zeros with multiplicity two which are given by $\Lambda_1^+\cup\Omega_2^-$ where the sets $\Lambda_1^+$ and $\Lambda_2^-$ are
\begin{equation}
\begin{array}{lll}
\Lambda_1^+=\{\mu_n\zeta_1\in l_{\zeta_1}\,(\mu_n>0,1\leq n\leq N<\infty\};\\
\Lambda_2^-=\{\nu_m\zeta_2\in l_{\zeta_2}\,(\nu_<0,1\leq m<M<\infty)\}.
\end{array}\label{eq2.17}
\end{equation}
\end{lemma}

P r o o f. If $w$ is a zero of the function $t_{0,0}(\lambda)$ \eqref{eq2.16}, then
\begin{equation}
u_0(w,0)(e_0^+(w,0))''-u'_0(w,0)(e_0^+(w,0))'+u''_0(w,0)\cdot e_0^+(w,0)=0.\label{eq2.18}
\end{equation}
Since
$$\frac d{d\lambda}t_{0,0}(\lambda)=\frac1{\lambda^4}\left\{\lambda^2\frac d{d\lambda}[u_0(\lambda,0)(e_0^+(\lambda,0))''-u'_0(\lambda,0)(e_0^+(\lambda,0))'+u''_0(\lambda,0)e_0^+(\lambda,0)]\right.$$
$$\left.-2\lambda[u_0(\lambda,0)(e_0^+(\lambda,0))''-u'_0(\lambda,0)(e_0^+(\lambda,0))'+u''_0(\lambda,0)e_0^+(\lambda,0)]\right\},$$
then, taking into account \eqref{eq2.18}, we obtain that
$$\frac d{d\lambda}t_{0,0}(\lambda)=\frac1{\lambda^2}\frac d{d\lambda}[u_0(\lambda,x)(e_0^+(\lambda,x))''-u'_0(\lambda,x)(e_0^+(\lambda,x))'+u''_0(\lambda,x)e_0^+(\lambda,x)].$$

Substitute in this equation $\lambda=w$, and taking into account that 
$$u'(w,0)(\stackrel{+}{v}(w,0))'=u(w,0)(\stackrel{+}{v}(w,0))''+u''(w,0)e^+(w,0)$$
in view of \eqref{eq2.18} we find that
$$\left.\frac d{d\lambda}t_{0,0}(\lambda)\right|_{\lambda=w}=\left.\frac1{w^2}\left\{u_0(w,0)(e_0^+(w,0))''\frac d{d\lambda}\ln\left(\frac{u(\lambda,0)(e^+(\lambda,0))''}{u'_0(\lambda,0)(e^+(\lambda,0)'}\right)\right\}\right|_{\lambda=w}$$
$$+u''_0(w,0)e^+(w,0)\left.\frac d{d\lambda}\ln\left(\frac{u''(\lambda,0)e^+(\lambda,0)}{u'(\lambda,0)(e^+(\lambda,0))'}\right)\right|_{\lambda=w}$$
$$=\frac{u_0(w,0)(e_0^+(w,0))}{w^2}\left.\frac d{d\lambda}\left\{\frac{u(\lambda,0)(e^+(\lambda,0))''+u''_0(\lambda,0)e_0^+(\lambda,0)}{u'_0(\lambda,0)(e^+(\lambda,0))'}\right\}\right|_{\lambda=w}.$$

\begin{remark}\label{r2.3}
The last formula implies that zero $\lambda=w$ of equation \eqref{eq2.18} is simple if only
$$u'_0(w,0)(e_0^+(w,0))'\left.\frac d{d\lambda}\left\{\frac{u(\lambda,0)(e^+(\lambda,0))''+u''_0(\lambda,0)e_0^+(\lambda,0)}{u'_0(\lambda,0)(e_0^+(\lambda,0))'}\right\}\right|_{\lambda=w}\not=0.$$
But if
$$u'_0(w,0)(e_0^+(w,0))'\left.\frac d{d\lambda}\left\{\frac{u(\lambda,0)e^+(\lambda,0)''+u''_0(\lambda,0)e_0^+(\lambda,0)}{u'(\lambda,0)(e_0^+(\lambda,0))'}\right\}\right|_{\lambda=w}=0,$$
then zero $w$ is of multiplicity $2$.

The case of simple zeros is studied in the works \cite{9,10}. We emphasize here zeros of multiplicity two of equation \eqref{eq2.18} and assume that all zeros of equation \eqref{eq2.18} are of multiplicity two.
\end{remark}

In order to check that the set of zeros of the function $t_{0,0}(\lambda)$ is finite, it is sufficient, due to analyticity of $t_{0,0}(\lambda)$ in sector $\Omega_0^-$, to ascertain that $t_{0,0}(\lambda)\rightarrow1$ as $\lambda\rightarrow1$ ($\lambda\in\Omega_0$). Really, every summand $u_0(\lambda,x)(e_0^+(\lambda,x))''$, $u'_0(\lambda,x)(e_0^+(\lambda,x))'$, and $u''_0(\lambda,x)e_0^+(\lambda,x)$, for $\lambda\rightarrow\infty$ ($\lambda\in\Omega_0$), tends to  $-\lambda^2+c(\lambda)$ where $c(\lambda)$ is a bounded function for all $\lambda\in\Omega_0^-$. So, $t_{0,0}(\lambda)$ \eqref{eq2.16} tends to 1 as $\lambda\rightarrow\infty$ ($\lambda\in\Omega_0^-$), which concludes the proof. $\blacksquare$

\begin{picture}(200,200)
\put(0,100){\vector(1,0){200}}
\put(100,0){\vector(0,1){200}}
\put(150,0){\vector(-1,2){100}}
\put(150,200){\vector(-1,-2){100}}
\put(200,67){\line(-3,1){200}}
\put(0,67){\line(3,1){200}}
%\put(30,190){$l_{\zeta_2}$}
\put(135,150){$\nu_m\zeta_2$}
%\put(80,170){$\Omega_2$}
\put(102,190){$\Omega_0^-$}
%\put(150,180){$\widehat{l}_{\zeta_3}$}
\put(123,150){$\circ$}
\put(70,147){$\times$}
\put(150,96){$\times$}
\put(154,105){$\mu_n\zeta_0$}
\put(138,47){$\nu_m\zeta_1$}
%\put(140,30){$\widehat{l}_{\zeta_2}$}
\put(125,40){$\circ$}
%\put(110,35){$\lambda'_n$}
\put(68,40){$\times$}
\put(40,33){$\mu_n\zeta_2$}
\put(60,5){$\Omega_1^-$}
%\put(60,60){$\Omega_3$}
%\put(10,85){$\nu_m\zeta_0$}
\put(20,105){$\nu_m\zeta_0$}
\put(50,96){$\circ$}
%\put(38,105){$\lambda''_n$}
%\put(190,105){$l_{\zeta_1}$}
\put(170,35){$\Omega_2^-$}
\put(45,150){$\mu_n\zeta_1$}
\qbezier(123,109)(100,130)(82,106)
\qbezier(100,80)(120,90)(120,107)
\qbezier(80,107)(80,80)(100,77)
\end{picture}

\hspace{20mm} Fig. 2

Geometric arrangement of zeros of the functions $t_{0,0}(\lambda)$, $t_{1,1}(\lambda)$, $t_{2,2}(\lambda)$ in the sectors $\Omega_0^-$, $\Omega_1^-$, $\Omega_2^-$ correspondingly is shown in Fig. 2. Besides, $\zeta_1\Lambda_1^+\cup\zeta_1\Lambda_2^-$ belongs to $\Omega_1^-$ and $\zeta_2\Lambda_1^+\cup\zeta_2\Lambda_2^-\subset\Omega_2^-$.
\vspace{5mm}

{\bf 2.4} Rewrite equality \eqref{eq1.4} as
\begin{equation}
r_0(\lambda)u_0(\lambda,x)=v_0(\lambda,x)+s_1(\lambda)v_1(\lambda,x)+s_2(\lambda)v_2(\lambda,x)\label{eq2.19}
\end{equation}
where
\begin{equation}
r_0(\lambda)=\frac1{t_{0,0}(\lambda)};\quad s_1(\lambda)\stackrel{\rm def}{=}\frac{t_{0,1}(\lambda)}{t_{0,0}(\lambda)};\quad s_2(\lambda)\stackrel{\rm def}{=}\frac{t_{0,2}(\lambda)}{t_{0,0}(\lambda)}.\label{eq2.19}
\end{equation}

\begin{remark}\label{r2.3}
Due to {\rm (a)} \eqref{eq1.13}, function $r_0(\lambda)u_0(\lambda,x)$ \eqref{eq2.19} behaves in the following way as $x\rightarrow\infty$:
$$r_0(\lambda)u_0(\lambda,x)\rightarrow e^{-i\lambda\zeta_0x}+s_1(\lambda)e^{-i\lambda\zeta_1x}+s_2(\lambda)e^{-i\lambda\zeta_2x}\quad(x\rightarrow\infty),$$
and thus the {\bf incident} (from $+\infty$) wave $e^{-i\lambda\zeta_0x}$ has the {\bf reflected (scattered)} wave $s_1(\lambda)e^{-i\lambda\zeta_1x}+s_2(\lambda)e^{-i\lambda\zeta_2x}$ where $s_1(\lambda)$ and $s_2(\lambda)$ are the {\bf scattering coefficients}. For $x\rightarrow-\infty$, function $r_0(\lambda)u_0(\lambda,x)$ \eqref{eq2.19}, due to {\rm (b)} \eqref{eq1.13}, has the asymptotis
$$r_0(\lambda)u_0(\lambda,x)\rightarrow r_0(\lambda)e^{-i\lambda\zeta_0x}\quad(x\rightarrow-\infty),$$
therefore, it is natural to consider $r_0(\lambda)$ as a {\bf transition coefficient} of the wave $e^{-i\lambda\zeta_0x}$.

Upon equating the (1,1)-elements in both sides of the equality {\rm (b)} \eqref{eq2.15}, we obtain
$$\zeta_p=t_{0,0}(\lambda)t_{0,0}^+(\lambda)+\zeta_1t_{0,1}t_{0,1}^+(\lambda)+\zeta_2t_{0,2}t_{0,2}^+(\lambda),$$
and thus
\begin{equation}
\zeta_pr_0(\lambda)r_0^+(\lambda)=1+\zeta_1s_1(\lambda)s_1^+(\lambda)+\zeta_2s_2(\lambda)s_2^+(\lambda).\label{eq2.21}
\end{equation}
One ought to consider this equality as a {\bf unitarity property of the scattering problem} for the incident wave $e^{-i\lambda\zeta_0x}$.
\end{remark}

Equation \eqref{eq2.19} implies two more scattering problems for the incident (from ``$+\infty$'') waves $e^{-i\zeta_1\lambda x}$ and $e^{-i\lambda\zeta_2x}$:
\begin{equation}
\begin{array}{ccc}
({\rm i})\quad r_1(\lambda)u_1(\lambda,x)=s_2(\lambda\zeta_1)v_0(\lambda,x)+v_1(\lambda,x)+s_1(\lambda\zeta_1)v_2(\lambda,x);\\
({\rm ii})\quad r_2(\lambda)u_2(\lambda,x)=s_1(\lambda\zeta_2)v_0(\lambda,x)+s_2(\lambda\zeta_2)v_1(\lambda,x)+v_2(\lambda,x),
\end{array}\label{eq2.22}
\end{equation}
besides, $r_k(\lambda)=r_0(\zeta_k\lambda)$ due to \eqref{eq2.5}.

Define the matrix functions
\begin{equation}
R(\lambda)\stackrel{\rm def}{=}\left[
\begin{array}{ccc}
r_0(\lambda)&0&0\\
0&r_1(\lambda)&0\\
0&0&r_2(\lambda)
\end{array}\right];\quad S(\lambda)\stackrel{\rm def}{=}\left[
\begin{array}{ccc}
1&s_1(\lambda)&s_2(\lambda)\\
s_2(\lambda\zeta_1)&1&s_1(\lambda\zeta_1)\\
-s_1(\lambda\zeta_2)&s_2(\lambda\zeta_2)&1
\end{array}\right],\label{eq2.23}
\end{equation}
then it is obvious that
$$T(\lambda)=R^{-1}(\lambda)S(\lambda).$$
Using (b) \eqref{eq2.15}, we obtain
\begin{equation}
\zeta_pR(\lambda)JR^*(\overline{\lambda})=S(\lambda)JS^*(\overline{\lambda}).\label{eq2.24}
\end{equation}
This equality one ought to consider as a {\bf unitarity property of the scattering problem} ({\bf conservation law of energetic balance}), besides, $R(\lambda)$ is the {\bf passage matrix} and $S(\lambda)$ is the {\bf scattering matrix}.
\vspace{5mm}

{\bf 2.4} Using \eqref{eq2.19}, calculate Wronskian $W\{u_0(\lambda,x),v_k(\lambda,x)\}$ ($k=1$, 2), then, taking \eqref{eq2.8} into account, we have
\begin{equation}
\begin{array}{lll}
({\rm i})\quad f_{0,1}(\lambda,x)=\sqrt3\lambda\zeta_2v_2^+(\lambda,x)-s_2(\lambda)\sqrt3\lambda\zeta_0v_0^+(\lambda,x)\quad(i\widehat{l}_{\zeta_2});\\
({\rm ii})\quad f_{0,2}(\lambda,x)=-\sqrt3\lambda\zeta_1v_1^+(\lambda,x)-s_1(\lambda)\sqrt3\lambda\zeta_0v_0^+(\lambda,x)\quad(i\widehat{l}_{\zeta_1})
\end{array}\label{eq2.25}
\end{equation}
where
\begin{equation}
f_{0,k}(\lambda,x)=r_0(\lambda)W\{u_0(\lambda,x),v_k(\lambda,x)\}\quad(k=1,2).\label{eq2.26}
\end{equation}

\begin{picture}(200,200)
\put(0,100){\vector(1,0){200}}
\put(100,0){\vector(0,1){200}}
\put(0,133){\vector(3,-1){200}}
\put(200,133){\vector(-3,-1){200}}
%\put(105,185){$il_{\zeta_0}$}
\put(120,150){$f_{0,2}(\lambda,x)$}
\put(45,150){$f_{0,1}(\lambda,x)$}
\qbezier(153,120)(100,140)(49,120)
\qbezier(100,23)(155,35)(168,125)
\qbezier(45,120)(30,70)(100,25)
%\qbezier(125,100)(115,110)(120,104)
%\put(140,105){$S_4(i)$}
%\put(70,110){$S_0(i)$}
%\put(80,88){$S_4$}
%\put(70,70){$S_2(i)$}
%\put(50,100){$S_1(i)$}
\put(20,40){$v_2^+(\lambda,x)$}
%\put(0,50){$il_{\zeta_1}$}
%\put(115,70){$S_3(i)$}
%\put(105,110){$S_5(i)$}
\put(160,50){$v_1^+(\lambda,x)$}
\put(185,135){$i\widehat{l}_{\zeta_1}$}
\put(0,145){$i\widehat{l}_{\zeta_2}$}
%\put(165,50){$\Omega_1$}
%\put(80,40){$\psi_1^*(\lambda,0)$}
\end{picture}

\hspace{20mm} Fig. 3

Functions $v_1^+(\lambda,x)$ and $v_2^+(\lambda,x)$ are holomorphic in the sectors $\Omega_2^-$ and $\Omega_1^-$ correspondingly, and $f_{0,k}(\lambda,x)$ are analytic in the sectors $\Omega_0^-\cap\Omega_k$ ($k=1$, 2). As a result, we have two jump problems on the rays $i\widehat{l}_{\zeta_1}$ and $i\widehat{l}_{\zeta_2}$ (Fig. 3).

Analogously, using (i) \eqref{eq2.22}, calculate Wronskians $W\{u_1(\lambda,x),v_k(\lambda,x)\}$ ($k=0$, 2), then, taking into account \eqref{eq2.8}, we obtain
\begin{equation}
\begin{array}{lll}
({\rm i})\quad f_{1,0}(\lambda,x)=-\sqrt3\lambda\zeta_2v_2^+(\lambda,x)+s_1(\lambda\zeta_1)\sqrt3\lambda\zeta_1v_1^+(\lambda,x);\quad(i\widehat{l}_{\zeta_0})\\
({\rm ii})\quad f_{1,2}(\lambda,x)=\sqrt3\lambda\zeta_0v_0^+(\lambda,x)-s_2(\lambda\zeta_1)\sqrt3\lambda\zeta_1v_1^+(\lambda,x);\quad(i\widehat{l}_{\zeta_1})
\end{array}\label{eq2.27}
\end{equation}
where
\begin{equation}
f_{1,k}(\lambda,x)=r_1(\lambda)W\{u_1(\lambda,x),v_k(\lambda,x)\}\quad(k=0,2)\label{eq2.28}
\end{equation}

\begin{picture}(300,200)
\put(70,100){\vector(1,0){200}}
\put(170,0){\vector(0,1){200}}
%\put(170,100){\vector(1,2){50}}
\put(260,131){\vector(-3,-1){200}}
\put(80,131){\vector(3,-1){200}}
%\put(168,170){$-\widehat\varkappa\zeta_2$}
\put(172,10){$i\widehat{l}_{\zeta_0}$}
%\put(198,162){$\circ$}
%\put(210,155){$r_0(\lambda\zeta_2)f_{2,0}(\lambda,x)$}
\put(200,140){$v_0^+(\lambda,x)$}
\put(190,50){$f_{1,0}(\lambda,x)$}
\put(90,50){$v_2^+(\lambda,x)$}
%\put(148,83){$\Omega_1^-$}
%\put(148,103){$\Omega_1^-$}
%\put(72,57){$il_{\zeta_1}$}
%\put(224,104){$\varkappa\zeta_0$}
\put(260,135){$i\widehat{l}_{\zeta_1}$}
\put(220,107){$f_{1,2}(\lambda,x)$}
%\put(260,80){$\Omega_1\cap\Omega_2^-$}
%\put(220,96){$\times$}
%\put(104,125){$v_2^+(\lambda,x)$}
%\put(247,58){$il_{\zeta_2}$}
\qbezier(220,85)(206,60)(170,65)
\qbezier(195,110)(220,120)(202,91)
\qbezier(203,112)(170,150)(139,112)
\qbezier(135,115)(130,80)(170,70)
\end{picture}

\hspace{20mm} Fig. 4.

Functions $v_0^+(\lambda,x)$ and $v_2^+(\lambda,x)$ are holomorphic at the sectors $\Omega_0^-$ and $\Omega_1^-$ correspondingly, and functions $f_{1,k}(\lambda,x)$, at the sectors $\Omega_1^-\cap\Omega_k$ ($k=0$, 2). Here we again have two jump problems on the rays $i\widehat{l}_{\zeta_0}$ and $i\widehat{l}_{\zeta_2}$ (Fig. 4). Also, using (ii) \eqref{eq2.22}, we find Wronskians $W\{u_2(\lambda,x),v_k(\lambda,x)\}$ ($k=0$, 1) and taking \eqref{eq2.8} into account, we have
\begin{equation}
\begin{array}{lll}
({\rm i})\quad f_{2,0}(\lambda,x)=\sqrt3\lambda\zeta_1v_1^+(\lambda,x)-s_2(\lambda\zeta_2)\sqrt3\lambda\zeta_2v_2^+(\lambda,x);\quad(i\widehat{l}_{\zeta_0})\\
({\rm ii})\quad f_{2,1}(\lambda,x)=-\sqrt3\lambda\zeta_0v_0^+(\lambda,x)-s_1(\lambda\zeta_2)\sqrt3\lambda\zeta_2v_2^+(\lambda,x);\quad(i\widehat{l}_{\zeta_2})
\end{array}\label{eq2.29}
\end{equation}
where
\begin{equation}
f_{2,k}(\lambda,x)=r_2(\lambda)W\{u_2(\lambda,x),v_k(\lambda,x)\}\quad(k=0,1).\label{eq2.30}
\end{equation}

\begin{picture}(200,200)
\put(0,100){\vector(1,0){200}}
\put(100,0){\vector(0,1){200}}
\put(200,67){\line(-3,1){200}}
\put(0,67){\line(3,1){200}}
\put(120,150){$v_0^+(\lambda,x)$}
%\put(175,140){$i\widehat{l}_{\zeta_1}$}
%\put(175,110){$\varphi_{1,2}(\lambda,x)$}
%\put(45,150){$g_{0,1}(\lambda,x)$}
\qbezier(153,120)(100,140)(49,120)
\qbezier(100,23)(155,35)(168,80)
\qbezier(60,115)(40,90)(60,85)
\qbezier(150,120)(168,110)(168,80)
\qbezier(44,80)(40,60)(100,25)
%\qbezier(125,100)(115,110)(120,104)
%\put(140,105){$S_2(i)$}
%\put(70,110){$S_4(i)$}
%\put(80,88){$S_4$}
%\put(70,70){$S_6(i)$}
%\put(50,100){$S_5(i)$}
%\put(0,52){$\widehat{\varphi}_{0,1}(\lambda,x)$}
\put(0,105){$f_{2,1}(\lambda,x)$}
%\put(0,135){$i\widehat{l}_{\zeta_2}$}
\put(140,30){$v_1^+(\lambda,x)$}
\put(30,27){$f_{2,0}(\lambda,x)$}
%\put(115,70){$S_1(i)$}
%\put(105,110){$S_3(i)$}
%\put(180,110){$g_{1,2}(\lambda,x)$}
%\put(165,50){$\Omega_1$}
%\put(80,40){$\psi_1^*(\lambda,0)$}
\end{picture}

\hspace{20mm} Fig. 5

Consider functions holomorphic in the corresponding sectors (Fig. 6):
$$\psi_0^+(\lambda,x)=v_0^+(\lambda,x)e^{i\lambda\zeta_0x}\quad(\lambda\in\Omega_0^-);$$
$$\psi_{1,2}(\lambda,x)=\frac1{\sqrt3\lambda}f_{1,2}(\lambda)e^{i\lambda\zeta_0x}\quad(\lambda\in\Omega_1^-\cap\Omega_2);$$
\begin{equation}
\psi_{2,1}(\lambda,x)=\frac1{\sqrt3\lambda}f_{2,1}(\lambda,x)e^{i\lambda\zeta_0x}\quad(\lambda\in\Omega_2^-\cap\Omega_1);\label{eq2.31}
\end{equation}
$$\psi_{1,0}(\lambda,x)=\frac1{\sqrt3\lambda}f_{1,0}(\lambda,x)e^{i\lambda\zeta_1x}\quad(\lambda\in\Omega_1^-\cap\Omega_0);$$
$$\psi_{2,0}(\lambda,x)=\frac1{\sqrt3\lambda}f_{2,0}(\lambda,x)e^{i\lambda\zeta_2x}\quad(\lambda\in\Omega_2^-\cap\Omega_0).$$
Equations (ii) \eqref{eq2.27} and (ii) \eqref{eq2.29} imply two jump problems on the rays $i\widehat{l}_{\zeta_1}$ and $i\widehat{l}_{\zeta_2}$:
\begin{equation}
\begin{array}{lll}
\psi_0^+(\lambda,x)-\psi_{1,2}(\lambda,x)=p_1(\lambda,x)\psi_1^+(\lambda,x)\quad(i\widehat{l}_{\zeta_1});\\
\psi_0^+(\lambda,x)+\psi_{2,1}(\lambda,x)=p_2(\lambda,x)\psi_2^+(\lambda,x)\quad(i\widehat{l}_{\zeta_2})
\end{array}\label{eq2.32}
\end{equation}
where
\begin{equation}
p_1(\lambda,x)\stackrel{\rm def}{=}\zeta_1s_2(\lambda\zeta_1)e^{-i\lambda\zeta_1x};\quad p_2(\lambda,x)\stackrel{\rm def}{=}-\zeta_2s_1(\lambda\zeta_2)e^{-i\lambda\zeta_2x}.\label{eq2.33}
\end{equation}
Since
$$\left.f_{1,2}(\lambda,x)\right|_{\lambda=i\tau\zeta_2\in il_{\zeta_2}}=r_0(i\tau)f_{0,1}(i\tau,x)\quad(\tau\geq0);$$
$$\left.f_{2,0}(\lambda,x)\right|_{\lambda=i\tau\zeta_1\in il_{\zeta_1}}=r_0(i\tau)f_{0,1}(i\tau,x)\quad(\tau\geq0),$$
then the boundary values of the function $\psi_{1,2}(\lambda,x)$ on the ray $il_{\zeta_2}$ coincide with the boundary values of the function $\psi_{2,0}(\lambda,x)$ on the ray $il_{\zeta_1}$. Analogously, it is proved that
$$\left.\psi_{2,1}(\lambda,x)\right|_{\lambda\in il_{\zeta_1}}=\left.\psi_{1,0}(\lambda,x)\right|_{\lambda\in il_{\zeta_2}}.$$

\begin{picture}(200,200)
\put(0,100){\vector(1,0){200}}
\put(100,200){\vector(0,-1){200}}
\put(200,67){\line(-3,1){200}}
\put(0,67){\line(3,1){200}}
\put(120,150){$\psi_0^+(\lambda,x)$}
\put(180,140){$i\widehat{l}_{\zeta_1}$}
\put(170,110){$\psi_{1,2}(\lambda,x)$}
%\put(45,150){$g_{0,1}(\lambda,x)$}
\qbezier(153,120)(100,140)(49,120)
\qbezier(100,23)(155,35)(168,80)
\qbezier(60,115)(40,90)(60,85)
\qbezier(150,120)(168,110)(168,80)
\qbezier(44,80)(40,60)(100,25)
%\qbezier(125,100)(115,110)(120,104)
%\put(140,105){$S_2(i)$}
%\put(70,110){$S_4(i)$}
%\put(80,88){$S_4$}
%\put(70,70){$S_6(i)$}
%\put(50,100){$S_5(i)$}
\put(0,42){$\widehat{\psi}_{1,0}(\lambda,x)$}
\put(0,105){$\psi_{2,1}(\lambda,x)$}
\put(140,25){$\widehat{\psi}_{2,0}(\lambda,x)$}
\put(110,0){$i\widehat{l}_{\zeta_0}$}
\put(0,140){$i\widehat{l}_{\zeta_2}$}
%\put(115,70){$S_1(i)$}
%\put(105,110){$S_3(i)$}
%\put(180,110){$g_{1,2}(\lambda,x)$}
%\put(165,50){$\Omega_1$}
%\put(80,40){$\psi_1^*(\lambda,0)$}
\end{picture}

\hspace{20mm} Fig. 7

By $\widehat{\psi}_{1,0}(\lambda,x)$, we denote a holomorphic in the sector $\Omega_2^-\cap\Omega_0$ function obtained from $\psi_{1,0}(\lambda,x)$ upon substituting $\lambda\rightarrow\overline{\lambda}$ where $\overline{\lambda}$ is the point symmetric to $\lambda$ relative to $i\widehat{l}_{\zeta_0}$, $\widehat{\psi}_{1,0}=\psi_{1,0}(\widehat{\lambda},x)$. Similarly, using symmetry, we define function $\widehat{\psi}_{2,0}(\lambda)$ in the sector $\Omega_1^-\cap\Omega_0$. As a result, we have holomorphic in the sectors $\Omega_1^-$ and $\Omega_2^-$ functions.

Multiply equality (i) \eqref{eq2.27} by $e^{i\lambda\zeta_1x}$ and equality (i) \eqref{eq2.29}, by $e^{i\lambda\zeta_2x}$,and subtract, then we obtain
\begin{equation}
\widehat{\psi}_{1,0}(\lambda,x)-\widehat{\psi}_{2,0}(\lambda,x)=p_3(\lambda,x)\psi_1^+(\lambda,x)+p_4(\lambda,x)\psi_2^+(\lambda,x)\label{eq2.34}
\end{equation}
where
\begin{equation}
p_3(\lambda,x)\stackrel{\rm def}{=}e^{i\lambda\zeta_2x}\zeta_1(s_1(\lambda\zeta_1)e^{i\lambda\zeta_1x}-e^{i\lambda\zeta_2x});\quad p_4(\lambda,x)\stackrel{\rm def}{=}e^{i\lambda\zeta_1x}\zeta_2(e^{i\lambda\zeta_1x}-s_2(\lambda\zeta_2)e^{i\lambda\zeta_2x}).\label{eq2.35}
\end{equation}
Equality \eqref{eq2.34} is the jump problem on the ray $i\widehat{l}_{\zeta_0}$.

So, on the rays $i\widehat{l}_{\zeta_1}$, $i\widehat{l}_{\zeta_2}$, and $i\widehat{l}_{\zeta_0}$ (Fig. 7), we have three jump problems \eqref{eq2.32}, \eqref{eq2.34} with summary parameters $\{p_k(\lambda,x)\}_1^4$ \eqref{eq2.33}, \eqref{eq2.35} that depend only on scattering coefficients $s_1(\lambda)$, $s_2(\lambda)$.

\begin{remark}\label{r2.4}
Functions $\psi_0^+(\lambda,x)$, $\psi_{1,2}(\lambda,x)$, $\psi_{2,1}(\lambda,x)$, $\widehat{\psi}_{1,0}(\lambda,x)$, $\widehat{\psi}_{2,0}(\lambda,x)$ tend to $1$ when $\lambda\rightarrow\infty$ (inside of the corresponding sector). For $\psi_0^+(\lambda,x)$, it is obvious, see \eqref{eq1.37}. Prove this statement for $\psi_{1,2}(\lambda,x)$, e.g. Using asymptotic behavior of $u_1(\lambda,x)$, $u'_1(\lambda,x)$ and $v_2(\lambda,x)$, $v'_2(\lambda,x)$ as $\lambda\rightarrow\infty$, we obtain
$$\psi_{1,2}(\lambda,x)=\frac{e^{i\lambda\zeta_0x}}{\sqrt3\lambda}\{(-i\lambda\zeta_1+o(1))(1+o(1))e^{-i\lambda\zeta_0x}-(1+o(1))(-i\lambda\zeta_2+o(1))e^{-i\lambda\zeta_0x}\}$$
$$=\frac1{\sqrt3\lambda}i\lambda(\zeta_2-\zeta_1)+o(1)=1+o(1).$$
\end{remark}

As a result, we have the Riemann boundary value problem on the contour formed by the rays $i\widehat{l}_{\zeta_0}$, $i\widehat{l}_{\zeta_1}$, $i\widehat{l}_{\zeta_2}$ for the piecewise holomorphic function
\begin{equation}
F(\lambda,x)\stackrel{\rm def}{=}\left\{
\begin{array}{lllll}
\psi_0^+(\lambda,x)\quad(\lambda\in\Omega_0^-);\\
\psi_{2,1}(\lambda,x)\quad(\lambda\in\Omega_2^-\cap\Omega_1);\\
\psi_{1,2}(\lambda,x)\quad(\lambda\in\Omega_1^-\cap\Omega_2);\\
\widehat{\psi}_{1,0}(\lambda,x)\quad(\lambda\in\Omega_2^-\cap\Omega_0);\\
\widehat{\psi}_{2,1}(\lambda,0)\quad(\lambda\in\Omega_1^-\cap\Omega_0).
\end{array}\right.\label{eq2.36}
\end{equation}
As is known \cite{17,18}, such function $F(\lambda,x)$ \eqref{eq2.36} is unambiguously restored from its jumps $\{p_k(\lambda,x)\}_1^4$ \eqref{eq2.33}, \eqref{eq2.35} using a Cauchy type integral, naturally, taking into account corresponding poles and their multiplicities, and also asymptotics for $\lambda\rightarrow\infty$,
$$F(\lambda,x)=1+\sum\limits_n\frac{R_n(\zeta_0,x)}{(\lambda-\mu_n\zeta_0)^2}+\sum\limits_n\frac{R_n(\zeta_2,x)}{(\lambda-\mu_n\zeta_2)^2}+\sum\limits_m\frac{\widehat{R}_m(\zeta_0,x)}{(\lambda-\nu_m\zeta_0)^2}$$
\begin{equation}
+\sum\limits_m\frac{\widehat{R}_m(\zeta,x)}{(\lambda-\nu_m\zeta_1)^2}+\frac1{2\pi i}\int\limits_{i\widehat{l}_{\zeta_1}}p_1(\mu,x)\psi_1^+(\mu,x)\frac{d\mu}{\mu-\lambda}+\int\limits_{il_{\zeta_2}}\frac{p_2(\mu,x)\psi_2^+(\mu,x)}{\mu-\lambda}d\mu\label{eq2.37}
\end{equation}
$$+\frac1{2\pi i}\int\limits_{i\widehat{l}_{\zeta_0}}\frac{p_3(\mu,x)\psi_1^+(\mu,x)+p_4(\mu,x)\psi_2(\mu,x)}{\mu-x}d\mu.$$
Coefficients $R_n(\zeta_2,x)$ are expressed via the coefficients $R_n(\zeta_0,x)$ by the formula
$$R_n(\zeta_2,x)=\zeta_1R_n(\zeta_0)e^{i\mu_n(\zeta_1-\zeta_2)x}=\zeta_1R_n(\zeta_0)e^{-\mu_n\sqrt3x}$$
which follows from the fact that function $g_{1,2}(\lambda,x)$, upon substituting $\lambda\rightarrow\lambda\zeta_1$, becomes $g_{2,0}(\lambda,x)$. Analogously,
$$R_m(\zeta_1,x)=\zeta_2\widehat{R}_m(\zeta_0,x)e^{i\nu_m(\zeta_2-\zeta_1)x}=\zeta_2\widehat{R}_m(\zeta_0,x)e^{\nu_m\sqrt3x}.$$
for $\lambda\in\Omega_0^-$, function $F(\lambda,x)$ \eqref{eq2.37}, coincides with $\psi_0(\lambda,x)$, therefore, taking into account these equalities, we have
$$\psi_0^+(\lambda,x)=1+\sum\limits_nR_n(\zeta_0,x)\left[\frac1{(\lambda-\mu_n)^2}+\frac{\zeta_1e^{-\sqrt3\mu_nx}}{(\lambda-\mu_n\zeta_2)^2}\right]+\sum\limits_m\widehat{R}_m(\zeta_0,x)\left[\frac1{(\lambda-\nu_m)^2}\right.$$
\begin{equation}
\left.+\frac{\zeta_2e^{\sqrt3\nu_mx}}{(\lambda-\nu_n\zeta_1)^2}\right]+\frac1{2\pi i}\int\limits_0^\infty\frac{p_1(-i\tau\zeta_1,x)\psi_2^+(-i\tau,x)}{\tau-i\zeta_2\lambda}d\tau+\frac1{2\pi i}\int\limits_0^\infty\frac{p_2(-i\tau\zeta_2,x)\psi_1^+(-i\tau,x)}{\tau-i\zeta_1\lambda}d\tau\label{eq2.38}
\end{equation}
$$+\frac1{2\pi i}\int\limits_0^\infty\frac{p_3(-i\tau,x)\psi_1^+(-i\tau,x)+p_4(-i\tau,x)\psi_2(-i\tau,x)}{\tau-i\lambda}d\tau\quad(\lambda\in\Omega_0^-).$$
Calculating boundary values $\lambda\rightarrow-it\zeta_1\in i\widehat{l}_{\zeta_1}$ in both sides of the equality and using Sokhotski formulas \cite{17,18}, we obtain
$$({\rm i})\quad\psi_1^+(-it,x)=1+\sum\limits_nR_n(\zeta_0,x)\left[\frac1{(it\zeta_1+\mu_n)^2}+\frac{\zeta_1e^{-\sqrt3\mu_nx}}{(it\zeta_1+\mu_n\zeta_2)^2}\right]$$
$$+\sum\limits_m\widehat{R}_m(\zeta_0,x)\left[\frac1{(it\zeta_1+\nu_m)^2}+\frac{\zeta_2e^{\sqrt3\nu_mx}}{(it\zeta_1+\nu_m\zeta_1)^2}\right]$$
\begin{equation}
+\frac12p_1(-it\zeta_1,x)\psi_2(-it,x)+\frac1{2\pi i}/\hspace{-4.4mm}\int\limits_0^\infty\frac{p_1(-i\tau\zeta_1,x)\psi_2^+(-i\tau,x)}{\tau-t}d\tau\label{eq2.39}
\end{equation}
$$+\frac1{2\pi i}\int\limits_0^\infty\frac{p_2(-i\tau\zeta_2,x)\psi_1^+(-i\tau,x)}{\tau-\zeta_2t}$$
$$+\frac1{2\pi i}\int\limits_0^\infty\frac{p_3(-i\tau)\psi_1(-i\tau,x)+p_4(-i\tau,x)\psi_2(-i\tau,x)}{\tau-\zeta_1t}d\tau.$$
Analogous calculations of the boundary values $\lambda\rightarrow-it\zeta_2\in i\widehat{l}_{\zeta_2}$ give us
$$({\rm ii})\quad\psi_2^+(-it,x)=1+\sum\limits_nR_n(\zeta_0,x)\left[\frac1{(it\zeta_2+\mu_n)^2}+\frac{\zeta_2e^{\sqrt3\nu_mx}}{(it\zeta_2+\mu_n\zeta_2)^2}\right]+\sum\limits_m\widehat{R}_m(\zeta_0)$$
\begin{equation}
\begin{array}{ccc}
{\displaystyle\times\left[\frac1{(it\zeta_2+\nu_m)^2}+\frac{\zeta_2e^{\sqrt3\nu_mx}}{(it\zeta_2+\nu_m\zeta_1)^2}\right]+\frac1{2\pi i}\int\limits_0^\infty\frac{p_1(-i\tau\zeta_1,x)\psi_2(-i\tau,x)}{\tau-\zeta_1t}d\tau}\\
{\displaystyle+\frac12p_2(-it\zeta_2,x)\psi_1(-it,x)+\frac1{2\pi i}/\hspace{-4.4mm}\int\limits_0^\infty\frac{p_2(-i\tau\zeta_2,x)\psi_1(-i\tau,x)}{\tau-t}d\tau}
\end{array}\label{eq2.40}
\end{equation}
$$+\frac1{2\pi i}\int\limits_0^\infty\frac{p_3(-i\tau,x)\psi_1(-i\tau,x)+p_4(-i\tau,x)\psi_2(-i\tau,x)}{\tau-\zeta_2t}dt.$$
It is left to obtain $N+M$ more equations to define the coefficients $\{R_n(\zeta_0,x)\}_1^N$ and $\{\widehat{R}_m(\zeta_0,x)\}_1^M$. Multiply equality \eqref{eq2.38} by $(\lambda-\mu_p)^{-1}$ and integrate it along the circle with the center at the point $\mu_p$ of the radius $r\ll1$ which does not contain any other points, apart from $\mu_p$, then we obtain
$$({\rm iii})\quad0=1+\sum\limits_{n\not=p}R_n(\zeta_0,x)\frac1{(\mu_p-\mu_n)^2}+\sum\limits_nR_n(\zeta_0,x)\frac{\zeta_1e^{-\sqrt3\mu_nx}}{(\mu_p-\mu_n\zeta_2)^2}+\sum\widehat{R}_m(\zeta_0,x)$$
\begin{equation}
\times\left[\frac1{(\mu_p+\nu_m)^2}+\frac{\zeta_2e^{\sqrt3\nu_mx}}{(\mu_p+\nu_m\zeta_1)^2}\right]+\frac1{2\pi i}\int\limits_0^\infty\frac{p_1(-i\tau\zeta_1,x)\psi_2^+(-i\tau,x)}{\tau-i\zeta_2\mu_p}d\tau\label{eq2.41}
\end{equation}
$$+\frac1{2\pi i}\int\limits_0^\infty\frac{p_2(-i\tau\zeta_2,x)\psi_1^+(-i\tau,x)}{\tau-i\zeta_1\mu_p}d\tau$$
$$+\frac1{2\pi i}\int\limits_0^\infty\frac{p_3(-i\tau,x)\psi_1(-i\tau,x)+p_4(-i\tau,x)\psi_2(-i\tau,x)}{\tau-i\mu_p}dt$$
($1\leq p\leq N$).

Analogously, we obtain $M$ more equations,
$$({\rm iv})\,0=1+\sum\limits_nR_n(\zeta_0,x)\left[\frac1{(\nu_q-\mu_n)^2}+\frac{\zeta_1e^{-\sqrt3\mu_nx}}{(\nu_q-\mu_n\zeta_2)^2}\right]+\sum\limits_{m\not=q}\widehat{R}_m(\zeta_0,x)\left[\frac1{(\nu_q-\nu_m)^2}\right.$$
\begin{equation}
\left.+\frac{\zeta_2e^{\sqrt3\nu_mx}}{(\nu_q-\nu_m\zeta_1)^2}\right]+\frac1{2\pi i}\int\limits_0^\infty\frac{p_1(-i\tau\zeta_1,x)\psi_2(-i\tau,x)}{\tau-i\zeta_2\nu_q}+\frac1{2\pi i}\int\limits_0^\infty\frac{p_2(-i\tau\zeta_2,x)\psi_1^+(-i\tau,x)}{\tau-i\zeta_1\nu_q}d\tau\label{eq2.42}
\end{equation}
$$+\frac1{2\pi i}\int\limits_0^\infty\frac{p_3(-i\tau,x)\psi_1(-i\tau,x)+p_4(-i\tau,x)\psi_2(-i\tau,x)}{\tau-i\nu_q}d\tau\quad(1\leq q\leq M).$$

\begin{conclusion}
We obtained the {\bf closed system of singular integral equations} {\rm (i) -- (iv)} \eqref{eq2.39} -- \eqref{eq2.42} relative to the unknowns $\{\psi_1^+(\lambda,x)\}$, $\{\psi_2^+(\lambda,x)\}_1^2$ and $\{R_n(\zeta_0,x)\}_1^N$, $\{\widehat{R}_m(\zeta_0,x)\}$, with free parameters $\{p_k(\lambda,x)\}_1^4$ and $\{\mu_n\}_1^N$, $\{\nu_m\}_1^M$. This system is analogous to the well-known Marchenko equation for Sturm -- Liouville operators.
\end{conclusion}

\begin{conclusion}
Knowing solution to the system of equations \eqref{eq2.39} -- \eqref{eq2.42}, define function $\psi_0^+(\lambda,x)$ using formula \eqref{eq2.38} using which, due to {\rm (i)}, {\rm(ii)} \eqref{eq1.48} and \eqref{eq1.49}, potentials $p(x)$, $q(x)$ are recovered on the right half-axis ($x\in\mathbb{R}_+$).
\end{conclusion}
\vspace{5mm}

{\bf 2.4} In the conclusion of this section, calculate reflectionless potentials assuming that $s_1(\lambda)=s_2(\lambda)=0$ and $n=m=1$. Then (iii) \eqref{eq2.41} and (iv) \eqref{eq2.42} imply the following system of linear equations:
\begin{equation}
\left\{
\begin{array}{lll}
{\displaystyle-1=R_1(\zeta_0,x)\frac{e^{-\sqrt3\mu_1x}}{\zeta_1\mu_1}+\widehat{R}_1(\zeta_0,x)\left[\frac1{(\mu_1+\nu_1)^2}+\frac{\zeta_2}{(\mu_1+\nu_1\zeta_2)^2}e^{\sqrt3\nu_1x}\right];}\\
{\displaystyle-1=R_1(\zeta_0,x)\left[\frac1{(\mu_1-\nu_1)^2}+\frac{\zeta_1e^{-\sqrt3\mu_1x}}{(\nu_1-\mu_1\zeta_2)^2}\right]+\widehat{R}_1(\zeta_0,x)\frac{e^{\sqrt3\nu_1x}}{\nu_1\zeta_2}}
\end{array}\right.\label{eq2.49}
\end{equation}
relative to $R_1(\zeta_0,x)$ and $\widehat{R}_1(\zeta_0,x)$. Determinant of this system is
\begin{equation}
\begin{array}{ccc}
{\displaystyle\Delta(x)=\frac{e^{\sqrt3(\nu_1-\mu_1)x}}{\mu_1\nu_1}-\left[\frac1{(\mu_1+\nu_1)^2}+\frac{\zeta_2}{(\mu_1+\nu_1\zeta_2)^2}e^{\sqrt3\nu_1x}\right]}\\
{\displaystyle\times\left[\frac1{(\mu_1-\nu_1)^2}+\frac{\zeta_1e^{-\sqrt3\mu_1x}}{(\nu_1+\mu_1\zeta_2)^2}\right],}
\end{array}\label{eq2.50}
\end{equation}
and solution to system \eqref{eq2.49} is
\begin{equation}
\begin{array}{ccc}
{\displaystyle R_1(\zeta_0,x)=\frac1{\Delta(x)}\left\{\frac{e^{\sqrt3\nu_1x}}{\nu_1\zeta_2}-\frac1{\mu_1+\nu_1}-\frac{\zeta_2}{(\mu_1+\nu_1\zeta_2)^2}e^{\sqrt3\nu_1x}\right\};}\\
{\displaystyle\widehat{R}_1(\zeta_0,x)=\frac1{\Delta(x)}\left\{\frac{e^{-\sqrt3\mu_1x}}{\zeta_1\mu_1}-\frac1{(\mu_1-\nu_1)^2}-\frac{\zeta_1e^{-\sqrt3\mu_1x}}{(\nu_1-\mu_1\zeta_2)^2}\right\}}.
\end{array}\label{eq2.51}
\end{equation}
Therefore, function $\psi_0^+(\lambda,x)$ \eqref{eq2.38} equals
\begin{equation}
\begin{array}{ccc}
{\displaystyle\psi_0^+(\lambda,x)=1+R_1(\zeta_0,x)\left\{\frac1{(\lambda-\mu_1)^2}+\frac{\zeta_1e^{-\sqrt3\mu_1x}}{(\lambda-\mu_1\zeta_2)^2}\right\}+\widehat{R}_1(\zeta_0,x)\left[\frac1{(\lambda-\nu_1)}\right.}\\
{\displaystyle\left.+\frac{\zeta_2e^{\sqrt3\zeta_1x}}{\lambda-\nu_1\zeta_1}\right].}
\end{array}\label{eq2.53}
\end{equation}
Hence, using \eqref{eq1.48}, we find potentials $p(x)$ and $q(x)$ on the right half-axis $\mathbb{R}_+$.

\section{Dual scattering problem for the waves incident from ``$-\infty$''}

{\bf 3.1} Analogously to \eqref{eq2.1}, expand every Jost solution $\{v_k(\lambda,x)\}_0^2$ in other Jost solutions,
\begin{equation}
v_k(\lambda,x)=\sum\limits_l\widetilde{t}_{k,l}(\lambda)u_l(\lambda,x)\quad(0\leq k\leq2),\label{eq3.1}
\end{equation}
or
\begin{equation}
v(\lambda,x)=\widetilde{T}(\lambda)u(\lambda,x)\label{eq3.2}
\end{equation}
where $u(\lambda,x)\stackrel{\rm def}{=}\col[u_0(\lambda,x),u_1(\lambda,x),u_2(\lambda,x)]$; $v(\lambda,x)\stackrel{\rm def}{=}\col[v_0(\lambda,x),v_1(\lambda,x),$ $v_2(\lambda,x)]$, and $\widetilde{T}(\lambda)$ is the {\bf dual transition matrix}
\begin{equation}
\widetilde{T}(\lambda)=\left[
\begin{array}{ccc}
\widetilde{t}_{0,0}(\lambda)&\widetilde{t}_{0,1}(\lambda)&\widetilde{t}_{0,2}(\lambda)\\
\widetilde{t}_{1,0}(\lambda)&\widetilde{t}_{1,1}(\lambda)&\widetilde{t}_{1,2}(\lambda)\\
\widetilde{t}_{2,0}(\lambda)&\widetilde{t}_{2,1}(\lambda)&\widetilde{t}_{2,2}(\lambda)
\end{array}\right].\label{eq3.3}
\end{equation}

\begin{remark}\label{r3.1}
Equation \eqref{eq2.2} implies that $T(\lambda)$ and $\widetilde{T}(\lambda)$ are mutually reciprocal, $T(\lambda)\widetilde{T}(\lambda)=I$. Due to (a) \eqref{eq2.15}, the inverse to $T(\lambda)$ matrix $T^{-1}=\overline{\overline{\zeta}}_pJT^*(\overline{\lambda})J^*$, and thus
\begin{equation}
\widetilde{T}(\lambda)=\overline{\zeta}_p\left[
\begin{array}{ccc}
t_{0,0}^+(\lambda)&\zeta_2t_{1,0}^+(\lambda)&\zeta_1t_{2,0}^+(\lambda)\\
\zeta_1t_{0,1}^+(\lambda)&t_{2,2}^+(\lambda)&\zeta_2t_{2,1}^+(\lambda)\\
\zeta_2t_{0,2}^+(\lambda)&\zeta_1t_{1,2}^+(\lambda)&t_{2,2}^+(\lambda)
\end{array}\right],\label{eq3.4}
\end{equation}
and thus $\widetilde{t}_{k,s}(\lambda)$ is expressed via $t_{k,s}^+(\lambda)$. For the elements $\widetilde{t}_{k,s}(\lambda)$, relations analogous to \eqref{eq2.5} are true. Similarly to \eqref{eq2.8}, for Wronskians $W_{k,s}(u,\lambda,x)$ of the functions $u_k(\lambda,x)$ and $u_s(\lambda,x)$, the following equalities hold:
\begin{equation}
\begin{array}{ccc}
W_{1,2}(u,\lambda,x)=\sqrt3\lambda u_0^+(\lambda,x);\quad W_{2,0}(u,\lambda,x)=\sqrt3\zeta_1u_1^+(\lambda,x);\\
W_{0,1}(u,\lambda,x)=\sqrt3\lambda\zeta_2u_2^+(\lambda,x).
\end{array}\label{eq3.5}
\end{equation}
\end{remark}

Consider equality \eqref{eq3.1} for $k=0$,
\begin{equation}
v_0(\lambda,x)=\widehat{t}_{0,0}(\lambda)u_0(\lambda,x)+\widetilde{t}_{0,1}(\lambda)u_1(\lambda,x)+\widetilde{t}_{0,2}(\lambda)u_2(\lambda,x),\label{eq3.6}
\end{equation}
and other equalities in \eqref{eq3.1} follow from this one upon substituting $\lambda\rightarrow\lambda\zeta_1$, $\lambda\rightarrow\lambda\zeta_2$. Analogously to considerations of Section 2, rewrite equality \eqref{eq3.6} as
\begin{equation}
\widetilde{r}_0(\lambda)v_0(\lambda,x)=u_0(\lambda,x)+\widetilde{s}_1(\lambda)u_1(\lambda,x)+\widetilde{s}_2(\lambda)u_2(\lambda,x)\label{eq3.7}
\end{equation}
where
$$\widetilde{r}_0(\lambda)\stackrel{\rm def}{=}\frac1{\widetilde{t}_{0,0}(\lambda)};\quad\widetilde{s}_1(\lambda)\stackrel{\rm def}{=}\frac{\widetilde{t}_{0,1}(\lambda)}{\widetilde{t}_{0,0}(\lambda)};\quad\widetilde{s}_2(\lambda)\stackrel{\rm def}{=}\frac{\widetilde{t}_{0,2}(\lambda)}{\widetilde{t}_{0,0}(\lambda)}.$$

For $x\rightarrow-\infty$, according to the asymptotic (b) \eqref{eq1.13}, function $\widetilde{r}_0(\lambda)v_0(\lambda,x)$ behaves in the following way at $-\infty$:
$$\widetilde{r}_0(\lambda)\widetilde{v}_0(\lambda,x)\rightarrow e^{-i\lambda\zeta_0x}+\widetilde{s}_1(\lambda)e^{-i\lambda\zeta_1x}+\widetilde{s}_2(\lambda)e^{-i\lambda\zeta_2x}\quad(x\rightarrow-\infty),$$
and the {\bf incident} (from $-\infty$) {\bf wave} $e^{-i\lambda\zeta_0x}$ has the {\bf reflected (scattered) wave} $\widetilde{s}_1(\lambda)e^{-i\lambda\zeta_1x}+\widetilde{s}_2(\lambda)e^{-i\lambda\zeta_2x}$ where $\widetilde{s}_1(\lambda)$ and $\widetilde{s}_2(\lambda)$ are the {\bf scattering coefficients} of this wave. At ``$+\infty$'', due to (a) \eqref{eq1.13}, function $r_0(\lambda)v_0(\lambda,x)$ \eqref{eq3.7} has the following asymptotic:
$$\widetilde{r}_0(\lambda)v_0(\lambda,x)\rightarrow\widetilde{r}_0(\lambda)e^{-i\lambda\zeta_0x}\quad(x\rightarrow\infty),$$
and thus it is natural to consider $\widetilde{r}_0(\lambda)$ as the {\bf transmission coefficient} of wave $e^{i\lambda\zeta_0x}$. Since $\widetilde{T}(\lambda)T(\lambda)=I$ and $\widetilde{T}(\lambda)=\overline{\zeta_p^*}JT^*(\overline{\lambda})J^*$, then $T(\lambda)=\overline{\zeta}_pJ\widetilde{T}^*(\overline{\lambda})J^*$, then from equality $\widetilde{T}(\lambda)J\widetilde{T}^*(\overline{\lambda})J^*=\zeta_p$ we obtain
$$\zeta_p=\widetilde{t}_{0,0}(\lambda)\widetilde{t}_{0,0}^+(\lambda)+\zeta_1t_{0,1}(\lambda)t_{0,1}^+(\lambda)+\zeta_2\widetilde{t}_{0,2}\widetilde{t}_{0,2}^+(\lambda),$$
and thus
\begin{equation}
\zeta_pr_0(\lambda)\widetilde{r}_0^+(\lambda)=1+\zeta_1\widetilde{s}_1(\lambda)s_1^+(\lambda)+\zeta_2\widetilde{s}_2(\lambda)s_2^*(\lambda).\label{eq3.8}
\end{equation}
Equality \eqref{eq3.8} describes the unitarity property of the scattering problem for an incident (from $-\infty$) wave $e^{-i\lambda\zeta_0x}$.

Note that equality \eqref{eq3.7}, upon substituting $\lambda\rightarrow\lambda\zeta_1$, $\lambda\rightarrow\lambda\zeta_2$, implies the equalities
\begin{equation}
\begin{array}{ccc}
{\rm(i)}\quad\widetilde{r}_1(\lambda)v_1(\lambda,x)=s_2(\lambda\zeta_1)u_0(\lambda,x)+u_1(\lambda,x)+s_1(\lambda\zeta_1)u_2(\lambda,x);\\
{\rm(ii)}\quad r_2(\lambda)v_2(\lambda,x)=s_1(\lambda\zeta_2)u_0(\lambda,x)+s_2(\lambda\zeta_2)u_1(\lambda,x)+u_2(\lambda,x)
\end{array}\label{eq3.9}
\end{equation}
describing scattering of the incident (from $-\infty$) waves $e^{-i\lambda\zeta_1x}$ and $e^{-i\lambda\zeta_2x}$ correspondingly. For scattering problems \eqref{eq3.9}, analogously to \eqref{eq3.8}, conservation laws hold.
\vspace{5mm}

{\bf 3.2} Using \eqref{eq3.7}, calculate Wronskians $W_{0,k}\{v_0(\lambda,x),u_k(\lambda,x)\}$ ($k=1$, 2), then, taking into account \eqref{eq3.5}, we have
\begin{equation}
\begin{array}{lll}
{\rm(i)}\quad g_{0,1}(\lambda,x)=\sqrt3\lambda\zeta_2u_2^+(\lambda,x)-\widetilde{s}_2(\lambda)\sqrt3\lambda\zeta_0u_0^+(\lambda,x);\\
{\rm(ii)}\quad g_{0,2}(\lambda,x)=-\sqrt3\lambda\zeta_1u_1^+(\lambda,x)+\widetilde{s}_1(\lambda)\sqrt3\lambda\zeta_0u_0^+(\lambda,x)
\end{array}\label{eq3.10}
\end{equation}
where
\begin{equation}
g_{0,k}(\lambda,x)\stackrel{\rm def}{=}\widetilde{r}_0(\lambda)W\{v_0(\lambda,x),u_k(\lambda,x)\}\quad(k=1,2).\label{eq3.11}
\end{equation}

\begin{picture}(400,200)
\put(24,100){\vector(1,0){200}}
\put(124,0){\vector(0,1){200}}
%\put(124,100){\vector(1,-2){50}}
%\put(124,100){\vector(-1,-2){50}}
\put(224,67){\line(-3,1){200}}
\put(24,67){\line(3,1){200}}
\put(26,50){$il_{\zeta_1}$}
%\put(89,40){$\times$}
%\put(99,44){$\varkappa\zeta_2$}
%\put(144,50){$\circ$}
%\put(150,50){$-\widehat{\varkappa}\zeta_2$}
\put(150,40){$g_{0,2}(\lambda,x)$}
%\put(90,3){$l_{\zeta_2}$}
%\put(181,15){$r_0(\lambda)f_{0,1}(\lambda,x)$}
\put(214,70){$il_{\zeta_2}$}
\qbezier(124,60)(159,50)(184,80)
\qbezier(124,63)(89,50)(64,80)
\qbezier(124,170)(74,140)(61,78)
\qbezier(124,170)(174,140)(187,78)
%\put(88,130){$\Omega_1^-$}
%\put(134,130){$\Omega_2^-$}
\put(174,133){$u_1^+(\lambda,x)$}
\put(30,135){$u_2^+(\lambda,x)$}
\put(55,50){$g_{0,1}(\lambda,x)$}
%\put(0,20){$r_0(\lambda)f_{0,2}(\lambda,x)$}
\end{picture}

\hspace{20mm} Fig. 8

Functions $u_1^+(\lambda,x)$ and $u_2^+(\lambda,x)$ are holomorphic in the sectors $\Omega_1$ and $\Omega_2$, and the functions $g_{0,k}(\lambda,x)$, in $\Omega_0\cap\Omega_k^-$ ($k=1$, 2) correspondingly. And we have two jump problems on the rays $il_{\zeta_1}$ and $il_{\zeta_2}$. Similarly, using (i) \eqref{eq3.9}, find Wronskians $W\{v_1(\lambda,x),u_k(\lambda,x)\}$ ($k=0$, 2), then, according to \eqref{eq3.5}, we have
\begin{equation}
\begin{array}{lll}
{\rm(i)}\quad g_{1,0}(\lambda,x)=-\sqrt3\lambda\zeta_2u_2^+(\lambda,x)+\widetilde{s}_1(\lambda\zeta_1)\sqrt3\lambda\zeta_1u_1^+(\lambda,x);\\
{\rm(ii)}\quad g_{1,2}(\lambda,x)=\sqrt3\zeta_0u_0^+(\lambda,x)-\widetilde{s}_2(\lambda\zeta_1)\sqrt3\lambda\zeta_1u_1^+(\lambda,x)
\end{array}\label{eq3.12}
\end{equation}
where
\begin{equation}
g_{1,k}(\lambda,x)=W\{v_1(\lambda,x),u_k(\lambda,x)\}\quad(k=0,2).\label{eq3.13}
\end{equation}

\begin{picture}(300,200)
\put(70,100){\vector(1,0){200}}
\put(170,0){\vector(0,1){200}}
%\put(170,100){\vector(1,2){50}}
\put(260,131){\vector(-3,-1){200}}
\put(170,100){\vector(3,-1){100}}
%\put(168,170){$-\widehat\varkappa\zeta_2$}
\put(172,190){$il_{\zeta_0}$}
%\put(198,162){$\circ$}
%\put(210,155){$r_0(\lambda\zeta_2)f_{2,0}(\lambda,x)$}
\put(200,140){$g_{2,0}(\lambda,x)$}
\put(190,50){$u_0^+(\lambda,x)$}
%\put(148,83){$\Omega_1^-$}
%\put(148,103){$\Omega_1^-$}
%\put(72,57){$il_{\zeta_1}$}
%\put(224,104){$\varkappa\zeta_0$}
\put(260,120){$g_{1,2}(\lambda,x)$}
%\put(260,80){$\Omega_1\cap\Omega_2^-$}
%\put(220,96){$\times$}
\put(104,125){$u_2^+(\lambda,x)$}
\put(247,58){$il_{\zeta_2}$}
\qbezier(193,110)(206,140)(170,133)
\qbezier(195,110)(220,120)(202,91)
\qbezier(203,88)(170,50)(139,88)
\qbezier(135,93)(130,120)(170,130)
\end{picture}

\hspace{20mm} Fig. 9

Functions $u_2^+(\lambda,x)$ and $u_0^+(\lambda,x)$ are analytic in the sectors $\Omega_2$ and $\Omega_0$, and functions $g_{1,k}(\lambda,x)$, in sectors $\Omega_1\cap\Omega_k^-$ ($k=0$, 1). Consequently, we again have two jump problems on the rays $il_{\zeta_0}$ and $il_{\zeta_1}$.

Taking into account (ii) \eqref{eq3.9}, calculate Wronskians $W\{v_2(\lambda,x),u_k(\lambda,x)\}$ ($k=0$, 1), then using \eqref{eq3.5} we have
\begin{equation}
\begin{array}{lll}
{\rm(i)}\quad g_{2,0}(\lambda,x)=\sqrt3\lambda\zeta_1u_1^+(\lambda,x)-\widetilde{s}_2(\lambda\zeta_2)\sqrt3\lambda\zeta_2u_2^+(\lambda,x);\\
{\rm(ii)}\quad g_{2,1}(\lambda,x)=-\sqrt3\lambda\zeta_0u_0^+(\lambda,x)+\widetilde{s}_1(\lambda\zeta_2)\sqrt3\lambda\zeta_2u_2^+(\lambda,x)
\end{array}\label{eq3.14}
\end{equation}
where
\begin{equation}
g_{2,k}(\lambda,x)=\widetilde{r}_2(\lambda)W\{v_2(\lambda,x),u_k(\lambda,x)\}\quad(k=0,1).\label{eq3.15}
\end{equation}

\begin{picture}(300,200)
\put(70,100){\vector(1,0){200}}
\put(170,0){\vector(0,1){200}}
%\put(170,100){\vector(-1,2){50}}
\put(70,133){\vector(3,-1){200}}
\put(170,100){\vector(-3,-1){100}}
%\put(100,190){$l_{\zeta_1}$}
\put(185,150){$u_1^+(\lambda,x)$}
\put(172,190){$il_{\zeta_0}$}
%\put(140,147){$\times$}
%\put(147,147){$\varkappa\zeta_1$}
\put(120,50){$u_0^+(\lambda,x)$}
\put(76,105){$g_{2,1}(\lambda,x)$}
\put(72,57){$il_{\zeta_1}$}
%\put(0,105){$r_0(\lambda\zeta_1)f_{1,2}(\lambda,x)$}
%\put(120,96){$\circ$}
%\put(104,105){$-\widehat\varkappa\zeta_0$}
\put(104,135){$g_{2,0}(\lambda,x)$}
%\put(247,58){$il_{\zeta_2}$}
%\put(30,150){$r_0(\lambda\zeta_1)f_{1,0}(\lambda,x)$}
\qbezier(193,89)(206,140)(170,133)
\qbezier(203,90)(170,50)(139,88)
\qbezier(149,108)(130,105)(139,88)
\qbezier(150,107)(150,120)(170,130)
%\put(177,110){$\Omega_2^-$}
%\put(174,80){$\Omega_0^-$}
%\put(50,130){$i\widehat{l}_{\zeta_2}$}
\end{picture}

\hspace{20mm} Fig. 10.

Functions $u_1^+(\lambda,x)$ and $u_0^+(\lambda,x)$ are holomorphic in the sectors $\Omega_1$ and $\Omega_0$, and functions $g_{2,k}(\lambda,x)$, correspondingly, in $\Omega_2\cap\Omega_k^-$ ($k=0$, 1). We also arrive at jump problems on the rays $il_{\zeta_0}$ and $il_{\zeta_1}$.

\begin{picture}(200,200)
\put(0,100){\line(1,0){200}}
\put(100,0){\vector(0,1){200}}
\put(0,133){\vector(3,-1){200}}
\put(200,133){\vector(-3,-1){200}}
%\put(2,130){$il_{\zeta_2}$}
\put(105,150){$g_{1,0}(\lambda,x)$}
%\put(100,185){$il_{\zeta_0}$}
%\qbezier(153,120)(100,140)(49,120)
%\qbezier(100,10)(170,20)(190,131)
%\qbezier(45,120)(30,70)(100,25)
\qbezier(100,140)(135,145)(141,115)
\qbezier(140,85)(155,105)(143,116)
\qbezier(100,60)(135,70)(144,85)
\qbezier(100,58)(70,60)(57,82)
\qbezier(56,86)(50,100)(58,115)
\qbezier(100,141)(65,150)(57,116)
%\put(10,55){$il_{\zeta_1}$}
\put(0,105){$g_{2,1}(\lambda,x)$}
\put(20,145){$g_{2,0}(\lambda,x)$}
%\put(180,135){$(-il_{\zeta_2})$}
\put(170,110){$g_{1,2}(\lambda,x)$}
%\put(173,53){$il_{\zeta_2}$}
\put(105,40){$u_0^+(\lambda,x)$}
%\put(45,40){$g_{2,1}(\lambda,x)$}
\end{picture}

\hspace{20mm} Fig. 11

These considerations imply that we have holomorphic functions in the sectors (see Fig. 11), besides, on each of the rays $il_{\zeta_1}$, $i\widehat{l}_{\zeta_1}$, $il_{\zeta_2}$, $i\widehat{l}_{\zeta_2}$, and $il_{\zeta_0}$ we have a jump problem. Analogously to \eqref{eq2.31}, define functions holomorphic in the sectors:
\begin{equation}
\left[
\begin{array}{lllll}
\varphi_0^+(\lambda,x)=u_0^+(\lambda,x)e^{i\lambda\zeta_0x}\quad(\lambda\in\Omega_0);\\
{\displaystyle\varphi_{1,2}(\lambda,x)=\frac1{\sqrt3\lambda}g_{1,2}(\lambda,x)e^{i\lambda\zeta_0x}\quad(\lambda\in\Omega_1\cap\Omega_2^-);}\\
{\displaystyle\varphi_{2,1}(\lambda,x)=\frac1{\sqrt3\lambda}g_{2,1}(\lambda,x)e^{\i\lambda\zeta_0x}\quad(\lambda\in\Omega_2\cap\Omega_1^-);}\\
{\displaystyle\varphi_{1,0}(\lambda,x)=\frac1{\sqrt3\lambda}g_{1,0}(\lambda,x)e^{i\lambda\zeta_2x}\quad(\lambda\in\Omega_1\cap\Omega_0^-);}\\
{\displaystyle\varphi_{2,0}(\lambda,x)=\frac1{\sqrt3\lambda}f_{2,0}(\lambda,x)e^{i\lambda\zeta_1x}\quad(\lambda\in\Omega_2\cap\Omega_0^-).}
\end{array}\right.\label{eq3.16}
\end{equation}
As in Section 2, it is easy to show that
$$\left.\psi_{1,2}(\lambda,x)\right|_{\lambda\in i\widehat{l}_{\zeta_1}}=\left.f_{2,0}(\lambda,x)\right|_{\lambda\in i\widehat{l}_{\zeta_2}};\quad\left.g_{2,1}(\lambda,x)\right|_{\lambda\in i\widehat{l}_{\zeta_2}}=\left.g_{1,0}(\lambda,x)\right|_{\lambda\in i\widehat{l}_{\zeta_1}}$$

Therefore, defining function $\widehat{\varphi}_{2,0}(\lambda,x)$ in the sector $(\Omega_1\cap\Omega_0^-)$ by symmetry (see Section 2), and $\widehat{\varphi}_{1,0}(\lambda,x)$ in the sector $(\Omega_2\cap\Omega_0^-)$, we obtain the holomorphic in the sectors $\Omega_0$, $\Omega_1$ and $\Omega_2$ functions
\begin{equation}
\Phi(\lambda,x)\stackrel{\rm def}{=}\left\{
\begin{array}{lllll}
\varphi_0^+(\lambda,x)\quad(\lambda\in\Omega_0);\\
\varphi_{1,2}(\lambda,x)\quad(\lambda\in\Omega_2\cap\Omega_1^-);\\
\varphi_{2,1}(\lambda,x)\quad(\lambda\in\Omega_1\cap\Omega_2^-);\\
\widehat{\varphi}_{2,0}(\lambda,x)\quad(\lambda\in\Omega_1\cap\Omega_0^-);\\
\varphi_{1,0}(\lambda,x)\quad(\lambda\in\Omega_2\cap\Omega_0^-).
\end{array}\right.\label{eq3.17}
\end{equation}
Equations (ii) \eqref{eq3.12} and (ii) \eqref{eq3.14} imply two jump problems on the rays $il_{\zeta_1}$  and $il_{\zeta_2}$:
\begin{equation}
\begin{array}{lll}
\varphi_0^+(\lambda,x)-\varphi_{1,2}(\lambda,x)=\widetilde{p}_1(\lambda,x)\varphi_1^+(\lambda,x)\quad(il_{\zeta_2});\\
\varphi_0^+(\lambda,x)+\psi_{2,1}(\lambda,x)=\widetilde{p}_2(\lambda,x)\varphi_2^+(\lambda,x)\quad(il_{\zeta_1})
\end{array}\label{eq3.18}
\end{equation}
where
$$\widetilde{p}_1(\lambda,x)=\widetilde{s}_2(\lambda\zeta_1)\zeta_1e^{i\lambda(\zeta_0-\zeta_1)x};\quad p_2(\lambda,x)=\widetilde{s}_1(\lambda\zeta_2)\zeta_2e^{i\lambda(\zeta_0-\zeta_2)x}.$$
And to obtain the jump problem on the ray $il_{\zeta_0}$, multiply equality (i) \eqref{eq3.12} by $e^{i\lambda\zeta_1x}$, and (i) \eqref{eq3.14}, by $e^{i\lambda\zeta_2x}$, and subtract, then, as a result, we have
\begin{equation}
\varphi_{1,0}(\lambda,x)-\varphi_2(\lambda,x)=p_3(\lambda,x)\varphi_1^+(\lambda,x)+p_4(\lambda,x)\varphi_2^+(\lambda,x)\label{eq3.19}
\end{equation}
where
\begin{equation}
\begin{array}{lll}
p_3(\lambda,x)\stackrel{\rm def}{=}e^{i\lambda\zeta_1x}\left[e^{i\lambda\zeta_1x}-\zeta_1\widetilde{s}_1(\lambda\zeta_1)e^{i\lambda\zeta_2x}\right];\\
 p_4(\lambda,x)\stackrel{\rm def}{=}e^{i\lambda\zeta_2x}\left[s_2\widetilde{s}_2(\lambda\zeta_2)e^{i\lambda\zeta_1}-e^{i\lambda\zeta_2x}\right].
 \end{array}\label{eq3.20}
\end{equation}

\begin{remark}\label{r2.1}
Obviously, $\widetilde{t}_{0,0}(\lambda)=\widetilde{t}_{00}^+(\lambda)$, and thus zeros of $t_{00}(\lambda)$ are the same as of $t_{00}(\lambda)$. Consequently, poles of the function $\varphi_{1,2}(\lambda,x)$ are $\{\mu_n\zeta_0\}_1^N$, of the function $\varphi_{2,1}(\lambda,x)$, $\{\nu_m\zeta_0\}$, and of the functions $\widehat{\varphi}_{1,0}(\lambda,x)$ and $\varphi_{2,0}(\lambda,x)$, correspondingly, $\{\mu_k\zeta_1\}$ and $\{\nu_n\zeta_2\}$.

Also note that every function from \eqref{eq3.17} tends to $1$ inside of the corresponding sector as $\lambda\rightarrow\infty$.
\end{remark}

Define a piecewise analytic function
\begin{equation}
\Phi(\lambda,x)\stackrel{\rm def}{=}\left\{
\begin{array}{lllll}
\varphi_0^+(\lambda,x)\quad(\lambda\in\Omega_0);\\
\varphi_{1,2}(\lambda,x)\quad(\lambda\in\Omega_1\cap\Omega_2^-);\\
\varphi_{2,1}(\lambda,x)\quad(\lambda\in\Omega_2\cap\Omega_1^-);\\
\widehat{\varphi}_{2,0}(\lambda,x)\quad(\lambda\in\Omega_1\cap\Omega_0^-);\\
\varphi_{1,0}(\lambda,x)\quad(\lambda\in\Omega_2\cap\Omega_0^-)
\end{array}\right.\label{eq3.21}
\end{equation}
which is holomorphic in the sectors $\Omega_0$, $\Omega_1$, $\Omega_2$ and on the rays $il_{\zeta_0}$, $il_{\zeta_1}$, $il_{\zeta_2}$ satisfies jump problems \eqref{eq3.18}, \eqref{eq3.19}. As is noted above, \cite{17, 18}, function $\Phi(\lambda,x)$ is defined unambiguously by the jumps on the rays $il_{\zeta_0}$, $il_{\zeta_1}$, $il_{\zeta_2}$ by a Cauchy type integral, taking into account poles and their multiplicity, and also its own asymptotic behavior at infinity,
$$\Phi(\lambda,x)=1+\sum\limits_n\frac{R'_n(\zeta_0,x)}{(\lambda-\mu_n\zeta_0)^2}+\sum\limits_n\frac{R'_n(\zeta_1,x)}{(\lambda-\mu_n\zeta_1)^2}+\sum\limits_m\frac{\widehat{R}'_m(\zeta_0,x)}{(\lambda-\nu_m\zeta_0)^2}+\sum\limits_m
\frac{\widehat{R}'(\zeta_2,x)}{(\lambda-\nu_m\zeta_2)^2}$$
\begin{equation}
+\frac1{2\pi i}\int\limits_{il_{\zeta_2}}\frac{p_1(\mu,x)\varphi'_1(\mu,x)}{\mu-\lambda}d\mu+\frac1{2\pi i}\int\limits_{il_{\zeta_1}}\frac{p_2(\mu,x)\varphi_2^+(\mu,x)}{\mu-x}d\mu\label{eq3.22}
\end{equation}
$$+\frac1{2\pi i}\int\limits_{l_{\zeta_0}}\frac{p_3(\mu,x)\varphi_1^+(\mu,x)+p_4(\mu,x)\varphi_2^+(\mu,x)}{\mu-\lambda}d\mu.$$

Since $\left.g_{1,2}(\lambda,x)\right|_{\lambda=\lambda\zeta_1}=g_{2,0}(\lambda,x)$ and $\left.g_{2,1}(\lambda,x)\right|_{\lambda=\lambda\zeta_2}=g_{1,0}(\lambda,x)$, then, as in Section 2, it is easy to show that
$$R'_n(\zeta_1,x)=R'_n(\zeta_0)\zeta_2e^{i\mu_n(\zeta_2-\zeta_1)x};\quad\widehat{R}'_m(\zeta_2,x)=R'_m(\zeta_0,x)\zeta_1e^{i\nu_m(\zeta_1-\zeta_2)x}.$$
Using these relations and the fact that function $\Phi(\lambda,x)$ \eqref{eq3.22}, for $\lambda\in\Omega_0$, coincides with function $\varphi_0(\lambda,x)$, we have
$$\varphi_0^+(\lambda,x)=1+\sum\limits_nR'_n(\zeta_0,x)\left[\frac1{(\lambda-\mu_n)^2}+\frac{\zeta_2e^{\sqrt3\mu_nx}}{(\lambda-\mu_n\zeta_1)^2}\right]+\sum\limits_m\widehat{R}'_m(\zeta_0,x)\left[\frac1{(\lambda-\nu_m)^2}\right.$$
\begin{equation}
\left.+
\frac{\zeta_1e^{-\nu_m\sqrt3x}}{(\lambda-\nu_m\zeta_2)^2}\right]+\frac1{2\pi i}\int\limits_0^\infty\frac{p_1(i\tau\zeta_2,x)\varphi_2^+(i\tau,x)}{\tau+i\zeta_1\lambda}d\tau+\frac1{2\pi i}\int\limits_0^\infty\frac{p_2(i\tau\zeta_1,x)\varphi_1^+(i\tau,x)}{\tau+i\zeta_2\lambda}d\tau\label{eq3.23}
\end{equation}
$$+\frac1{2\pi i}\int\limits_0^\infty\frac{p_3(i\tau,x)\psi_1(i\tau,x)+p_4(i\tau,x)\varphi_2(i\tau,x)}{\tau+i\lambda}d\tau\quad(\lambda\in\Omega_0).$$
Calculating boundary values in both sides of formula \eqref{eq3.23} when $\lambda\rightarrow it\zeta_1\in il_{\zeta_1}$ and using Sokhotski formulas \cite{17, 18}, we obtain
$$({\rm i})\quad\varphi_2^+(\lambda,x)=1+\sum\limits_nR'_n(\zeta_0,x)\left[\frac1{(it\zeta_1-\mu_n)^2}+\frac{\zeta_2e^{\sqrt3\mu_nx}}{(it\zeta_1-\mu_n\zeta_1)^2}\right]$$
\begin{equation}
+\sum\limits_mR'(\zeta_0,x)\left[
\frac1{(it\zeta_1-\nu_m)^2}+\frac{\zeta_1e^{-\nu_m\sqrt3x}}{(it\zeta_1-\nu_m\zeta_2)^2}\right]+\frac1{2\pi i}\int\limits_0^\infty\frac{p_1(i\tau\zeta_2,x)\varphi_2^+(i\tau,x)}{\tau-\zeta_2t}d\tau\label{eq3.24}
\end{equation}
$$+\frac12(p_2(-it\zeta_2,x)\varphi_1^+(it,x))+\frac1{2\pi i}/\hspace{-4.4mm}\int\limits_0^\infty\frac{p_2(i\tau\zeta_1,x)\varphi_1^+(i\tau,x)}{\tau-t}d\tau$$
$$+\frac1{2\pi i}\int\limits_0^\infty\frac{p_3(i\tau,x)\varphi_1(i\tau,x)+p_4(i\tau,x)\varphi_2(i\tau,x)}{\tau-\zeta_1t}.$$
Analogously, boundary values $\lambda\rightarrow it\zeta_2\in il_{\zeta_2}$ in equality \eqref{eq3.23} give us
$$({\rm ii})\quad\varphi_1^+(\lambda,x)=1+\sum\limits_nR'_n(\zeta_0,x)\left[\frac1{(it\zeta_2-\mu_n)^2}+\frac{\zeta_2e^{\sqrt3\mu_nx}}{(it\zeta_2-\mu_n\zeta_1)^2}\right]$$
\begin{equation}
+\sum\limits_m\widehat{R}'_m(\zeta_0,x)\left[\frac1{(it\zeta_2-\nu_m)^2}+\frac{\zeta_1e^{-\nu_m\sqrt3x}}{(it\zeta_2-\nu_m\zeta_2)^2}\right]+\frac12[p_1(it\zeta_2,x)\varphi_2^+(it,x)]\label{eq3.25}
\end{equation}
$$+\frac1{2\pi i}/\hspace{-4.4mm}\int\limits_0^\infty\frac{p_1(i\tau\zeta_2,x)\varphi_2^+(i\tau,x)}{\tau-t}d\tau+\frac1{2\pi i}\int\limits_0^\infty\frac{p_2(i\tau\zeta_1,x)\varphi_1^+(i\tau,x)}{\tau-\zeta_1t}d\tau$$
$$+\frac1{2\pi i}\int\limits_0^\infty\frac{p_3(i\tau,x)\varphi_1^+(i\tau,x)+p_4(i\tau,x)\varphi_2^+(i\tau,x)}{\tau-\zeta_2t}d\tau.$$
And it is necessary to obtain $N+M$ more equations in order to find coefficients $\{R'_n(\zeta_0)\}_1^N$ and $\{\widehat{R}'_m(\zeta_0,x)\}_1^M$. Multiplying \eqref{eq3.23} by $(\lambda-\mu_p)^{-1}$ and integrating over a circle with the center at $\mu_p$ of the radius $r\ll1$ which does not contain other points, apart from $\mu_p$, we obtain
\begin{equation}
\begin{array}{ccc}
{\displaystyle({\rm iii})\quad0=1+\sum\limits_{n\not=p}R'_n(\zeta_0,x)\frac1{(\mu_p-\mu_n)^2}+\sum\limits_nR'_n(\zeta_0,x)\frac{\zeta_2e^{\sqrt3\mu_nx}}{(\mu_p-\mu_n\zeta_1)^2}+\sum\limits_m\widehat{R}'_m(\zeta,x)}\\
{\displaystyle\times\left[\frac1{(\mu_p-\nu_m)^2}+\frac{\zeta_1ee^{-\nu_mx}}{(\mu_p-\nu_m\zeta_2)^2}\right]+\frac1{2\pi i}\int\limits_0^\infty\frac{p_1(i\tau\zeta_2,x)\varphi_2^+(i\tau,x)}{\tau+i\zeta_1\mu_p}d\tau}
\end{array}\label{eq3.26}
\end{equation}
$$+\frac1{2\pi i}\int\limits_0^\infty\frac{p_2(i\tau\zeta_1,x)\varphi_1^+(i\tau,x)}{\tau+i\zeta_2\mu_p}d\tau+\frac1{2\pi i}\int\limits_0^\infty\frac{p_3(i\tau,x)\varphi_1^+(i\tau,x)+p_4(i\tau,x)\varphi_2^+(i\tau,x)}{\tau+i\mu_p}d\tau$$
($1\leq p\leq N$).

Analogously, we obtain $M$ more equations
$$({\rm iv})\quad0=1+\sum\limits_nR'_n(\zeta_0,x)\left[\frac1{(\nu_q-\mu_n)^2}+\frac{\zeta_2e^{\sqrt3\mu_nx}}{(\nu_q-\mu_n\zeta_1)^2}\right]+\sum\limits_{m\not=q}\widehat{R}'_m(\zeta_0,x)\frac1{(\nu_q-\nu_m)^2}$$
\begin{equation}
+\sum\limits_m\widehat{R}'_m(\zeta_0,x)\frac{\zeta_1e^{-\nu_m\sqrt3x}}{(\nu_q-\nu_m\zeta_2)^2}+\frac1{2\pi i}\int\limits_0^\infty\frac{p_1(i\tau\zeta_2,x)\varphi_2^+(i\tau,x)}{\tau+i\zeta_1\nu_q}d\tau\label{eq3.27}
\end{equation}
$$+\frac1{2\pi i}\int\limits_0^\infty\frac{p_2(i\tau\zeta_1,x)\varphi_1^+(i\tau,x)}{\tau+i\zeta_2\nu_q}d\tau+\frac1{2\pi i}\int\limits_0^\infty\frac{p_3(i\tau,x)\varphi_1^+(i\tau,x)+p_4(i\tau,x)\varphi_2^+(i\tau,x)}{\tau+i\nu_q}d\tau$$
($1\leq q\leq M$).

\begin{conclusion}
We found the {\bf second closed} system of singular integral equations {\rm (i) -- (iv)} \eqref{eq3.24} -- \eqref{eq3.27} for unknowns $\varphi_1^+(\lambda,x)$, $\varphi_2^+(\lambda,x)$ and $\{R'_n(\zeta_0,x)\}_1^N$, $\{R'_m(\zeta_0,x)\}_1^M$, summary parameters of which are functions $\{p_k(\lambda,x)\}_1^4$ that are expressed via the scattering coefficients $\widetilde{s}_1(\lambda)$ and $\widetilde{s}_2(\lambda)$ and points $\{\mu_n\}_1^N$ and $\{\nu_m\}_1^M$. One ought to consider this system of equations as an analogue of the well-known Marchenko equation.
\end{conclusion}

\begin{conclusion}
Knowing solution to the system of equations \eqref{eq3.24} -- \eqref{eq3.27}, we define function $\varphi_0^+(\lambda,x)$ by formula \eqref{eq3.23}. Next, using (i), (ii) \eqref{eq1.55} and \eqref{eq1.56} we define potentials $p(x)$ and $q(x)$ for $x\in\mathbb{R}_-$.
\end{conclusion}

In the conclusion of this section, note that, analogously to Subsec. 4.3, we can calculate reflectionless potentials $p(x)$ and $q(x)$ on the left half-axis $\mathbb{R}_-$.

\renewcommand{\refname}{ \rm 
\end{Large}
\end{document}